\documentclass[preprint,12pt,authoryear]{elsarticle}

\usepackage{lipsum}
\usepackage{mathtools}
\usepackage{mathrsfs}
\usepackage{bbm}
\usepackage{graphicx}
\usepackage{amsmath}
\usepackage{amsthm}
\usepackage{amssymb}
\usepackage{amsfonts}
\usepackage{epsfig}
\usepackage{graphicx,float}
\usepackage{color}
\usepackage{bbm}
\usepackage[table]{xcolor}
\usepackage[margin=0.75in]{geometry}
\usepackage{comment}
\usepackage{enumitem}
\usepackage{wrapfig}
\usepackage{caption}
\usepackage{subcaption}
\usepackage{setspace}
\usepackage{algorithm}
\usepackage{qtree}
\usepackage{csquotes}
\usepackage[noend]{algpseudocode}
\usepackage{mdframed}
\usepackage{wrapfig}
\makeatletter
\def\BState{\State\hskip-\ALG@thistlm}
\makeatother
\doublespace
\usepackage{natbib} \citeindextrue

\usepackage[latin1]{inputenc}
\usepackage{tikz}
\usetikzlibrary{shapes,arrows}

\DeclareMathOperator*{\argmin}{\arg\!\min}
\DeclareMathOperator*{\argmax}{\arg\!\max}
\newmdtheoremenv{theo}{Theorem}

\usepackage{lineno}

\journal{European Journal of Operations Research}

\begin{document}

\begin{frontmatter}


\title{Stochastic Optimization with Parametric Cost Function Approximations}
\author[RP3]{Raymond T. Perkins III}
\author[WBP]{Warren B. Powell}

\address[RP3]{raymondp@princeton.edu}
\address[WBP]{powell@princeton.edu}
\begin{abstract}
A widely used heuristic for solving stochastic optimization problems is to use a
deterministic rolling horizon procedure, which has been modified to handle uncertainty (e.g. buffer stocks, schedule slack). This approach has been criticized for its use of a deterministic approximation of a stochastic problem, which is the major motivation for stochastic programming.  We recast this debate by identifying both deterministic and stochastic approaches as policies for solving a stochastic base model, which may be a simulator or the real world.  Stochastic lookahead models (stochastic programming) require a range of approximations to keep the problem tractable.  By contrast, so-called deterministic models are actually parametrically modified cost function approximations which use parametric adjustments to the objective function and/or the constraints.  These parameters are then optimized in a stochastic base model which does not require making any of the types of simplifications required by stochastic programming.  We formalize this strategy and describe a gradient-based stochastic search strategy to optimize the parameters.
\end{abstract}
\begin{keyword}
Stochastic Optimization \sep Stochastic Programming \sep Decisions under uncertainty \sep Parametric Cost Function Approximation \sep Cost Function Approximation \sep Policy Search
\end{keyword}

\end{frontmatter}

\section{Introduction}
There has been a long history in industry of using deterministic optimization models to make decisions that are then implemented in a stochastic setting.  Energy companies use deterministic forecasts of wind, solar and loads to plan energy generation (\cite{wallace2003stochastic}); airlines use deterministic estimates of flight times to schedule aircraft and crews (\cite{lan2006planning}); and retailers use deterministic estimates of demands and travel times to plan inventories (\cite{harrison1999multi}).  These models have been widely criticized in the research community for not accounting for uncertainty, which often motivates the use of large-scale stochastic programming models which explicitly model uncertainty in future outcomes (\cite{mulvey1995robust} and \cite{birge2011introduction}). These large-scale stochastic programs have been applied to unit commitment (\cite{jin2011modeling}), hydroelectric planning (\cite{carpentier2015managing}), and transportation (\cite{lium2009study}).  These models use large scenario trees to approximate potential future events, but result in very large-scale optimization models that can be quite hard to solve in practice.

We make the case that these previous approaches ignore the true problem that is being solved, which is always stochastic.  The so-called ``deterministic models'' used in industry are almost always parametrically modified deterministic approximations, where the modifications are designed to handle uncertainty.  Both the ``deterministic models'' and the ``stochastic models'' (formulated using the framework of stochastic programming) are examples of lookahead policies to solve a stochastic optimization problem.  The stochastic optimization problem is to find the best policy which is typically tested using a simulator, but may be field tested in an online environment (the real world).

In this paper, we characterize these modified deterministic models as {\it parametric cost function approximations} which puts them into the same category as other parameterized policies that are well known in the research community working on policy search (\cite{ng2000pegasus}, \cite{peshkin2000learning}, \cite{hu2007evolutionary}, \cite{Robots}, and \cite{mannor2003cross}).  Our use of modified linear programs is new to the policy search literature, where ``policies'' are typically parametric models such as linear models (``affine policies''), structured nonlinear models (such as (s,S) policies for inventories) or neural networks. The process of designing the modifications (in this paper, these modifications always appear in the constraints) requires the same art as the design of any statistical model or parametric policy.  The heart of this paper is the design of gradient-based search algorithms, which are nontrivial in this setting.

This paper formalizes the idea, used for years in industry, that an effective way to solve complex stochastic optimization problems is to shift the modeling of the stochastics from a lookahead approximation, where even deterministic lookahead models can be hard to solve, to the stochastic base model, typically implemented as a simulator but which might also be the real world.  Tuning a model in a stochastic simulator makes it possible to handle arbitrarily complex dynamics, avoiding the many approximations (such as two-stage models, exogenous information that is independent of decisions) that are standard in stochastic programming.

Parametric cost function approximations also make it possible to exploit structural properties.  For example, it may be obvious that the way to handle uncertainty when planning energy generators in a unit commitment problem is to require extra reserves at all times of the day.  A stochastic programming model encourages this behavior, but the requirement for a manageable number of scenarios will produce the required reserve only when one of the scenarios requires it.  Imposing a reserve constraint (which is a kind of cost function approximation) allows us to impose this requirement at all times of the day, and to tune this requirement under very realistic conditions.

Designing a parametric cost function approximation closely parallels the design of any parametric statistical model, which is part art (creating the model) and part science (fitting the model). To illustrate the process of designing a parametric cost function approximation, we use the setting of a time-dependent stochastic inventory planning problem that arises in the context of energy storage, but could represent any inventory planning setting.  We assume we have access to rolling forecasts where forecast errors are based on careful modeling of actual and predicted values for energy loads, generation from renewable sources, and prices.  The combination of the time-dependent nature and the availability of rolling forecasts which are updated each time period make this problem a natural setting for lookahead models, where the challenge is how to handle uncertainty.  We have selected this problem since it is relatively small, simplifying the extensive computational work.  However, our methodology is scalable to any problem setting which is currently being solved using a deterministic model.

Most important, the parametric CFA opens up a fundamentally new approach for providing practical tools for solving high-dimensional, stochastic programming problems.  It provides an alternative to classical stochastic programming with its focus on optimizing a stochastic lookahead model which requires a variety of approximations to make it computationally tractable.  The parametric CFA makes it possible to incorporate problem structure, such as the recognition that robust solutions can be achieved using standard methods such as schedule slack and/or buffer stocks. The parametric CFA makes it possible to incorporate problem structure for handling uncertainty.  Some examples include:
\begin{itemize}
	\item Supply chains handle uncertainty by introducing buffer stocks.
	\item Airlines handle uncertainty due to weather and congestion by using schedule slack.
	\item FedEx plans for equipment problems by maintaining spare aircraft at different locations around the country.
	\item Hospitals can handle uncertainty in blood donations and the demand for blood by maintaining supplies of O-minus blood, which can be used by anyone.
	\item Grid operators handle uncertainty in generator failures, as well as uncertainty in energy from wind and solar, by requiring generating reserves.
\end{itemize}
Central to our approach is the ability to manage uncertainty by recognizing effective strategies for responding to unexpected events.  We would argue that this structure is apparent in many settings, especially in complex resource allocation problems.  At a minimum, we offer that our approach represents an interesting, and very practical, alternative to stochastic programming.

This paper makes the following contributions.  1) We introduce and develop the idea of parameterized cost function approximations as a tool for solving important classes of stochastic optimization problems, shifting the focus from solving complex, stochastic lookahead models to optimizing a stochastic base model.  This approach is computationally comparable to solving deterministic approximations, with the exception that the parametric modifications have to be optimized, typically in a simulator that avoids the many approximations made in stochastic lookahead models. 2) We derive the policy gradients for parameterized right-hand sides using the properties of the underlying linear program.  3) We illustrate different styles of parametric approximations using the context of a nonstationary energy inventory problem, and quantify the benefits over a basic deterministic lookahead without adjustments.

Our presentation is organized as follows. The modeling framework is given in section 2.  We then provide an overview of the different classes of policies in section 3. We defer the literature review until section 3 which allows us to put the literature in the framework of the different classes of policies.  Alternative designs for parametric CFAs are presented in section 4, along with the derivation of the gradient of the base model with respect to the policy parameters.  Section 5 presents a series of numerical results.

\section{Canonical Model}
 Sequential, stochastic decision problems require a richer notation than standard linear programs and deterministic problems. For the sake of notational consistency, we follow the canonical model in \cite{powell2011approximate} which breaks dynamic programs into five dimensions:

\begin{itemize}		
	\item The {\bf state variable}, $S_t$, is the minimally dimensioned function of history that, combined with the exogenous information process, contains the information needed to compute the cost function, the constraints, and the transition function, from time $t$ onward. We use the initial state, $S_0$, to represent any deterministic parameters, as well as probabilistic beliefs (if needed).
	\item A {\bf decision}, $x_t$, is an $n$-dimensional vector that must satisfy $x_t \in \mathcal{X}_t$, which is typically a set of linear constraints.  Decisions are determined by a decision function (policy) which we denote by $X^\pi_t(S_t)$, where $\pi$ carries the information that determines the structure and parameters that define the function.  
	\item The {\bf exogenous information}, $W_t$, describes the information that first becomes known at time $t$.  We let $\omega \in \Omega$ be a sample path of $W_1, \ldots, W_T$.  Let $\mathcal{F}$ be the sigma algebra on $\Omega$, and let $\mathcal{P}$ be a probability measure on $(\Omega,\mathcal{F})$, giving us a probability space $(\Omega,\mathcal{F},\mathcal{P})$.  Next let $\mathcal{F}_t = \sigma(W_1, \ldots, W_t)$ be the sigma-algebra generated by $W_1, \ldots, W_t, $ where $ (\mathcal{F}_t)_{t=1}^T$ forms a filtration. The information $W_t$ may depend on the state $S_t$ and/or the action $x_t$, which means it depends on the policy.  If this is the case, we write our probability space as $(\Omega^\pi,\mathcal{F}^\pi,\mathcal{P}^\pi)$, with the associated expectation operator $\mathbb{E}^\pi$.
	\item The {\bf transition function}, $S^M(\cdot)$,  explicitly describes the relationship between the state of the model at time $t$ and $t+1$,
\begin{equation}
	S_{t+1} = S^M(S_t, x_t, W_{t+1}).
\end{equation}
	\item The {\bf objective function} is used to evaluate the effectiveness of a policy or sequence of decisions. It  minimizes the expected sum of the costs $C(S_t,x_t)$ in each time period $t$ over a finite horizon, where we seek to find the policy that solves,
\begin{equation} \label{objective1}
	\min_{\pi \in \Pi} \: \: \mathbb{E}^{\pi} \bigg[ \sum^T_{t=0} C(S_t, X^{\pi}_t(S_t)) \: \bigg| \: S_0  \bigg]
\end{equation}
where $S_{t+1} = S^M(S_t,X^\pi_t(S_t),W_{t+1}).$ We use $\mathbb{E}^{\pi}(\cdot)$ since the exogenous variables in the model may be affected by the decisions generated by our policy. Therefore we express the expectation as dependent on the policy. Since stochastic problems incorporate uncertainty in the model a variety of risk measures can be used in replacement of expectation. Equation \eqref{objective1}, along with the transition function and the exogenous information process, is called the {\it base model}.
\end{itemize}

This canonical model can be used to model virtually {\it any} sequential, stochastic decision problem as long as we are using expectations instead of risk measures.  We use this setting to put different policies onto a standard footing for comparison. In the next section we describe the major classes of policies that we can draw from to solve the problem.  We use this framework to review the literature.

\section{Solution Strategies}
The wide breadth of problems under the umbrella of \emph{optimization under uncertainty} has motivated an equally diverse assortment of policies and methods. This collection of strategies can be partitioned into two broad classifications: policy search and policies based on lookahead approximations. Approximate lookahead methods approximate the impact of a current decision on the future. This class  can be divided into two sub-categories: Approximations of lookahead models and value function approximations. Approximations of lookahead models  approximate the true model using a range of strategies. Value functions (typically known as value function approximations) use some form of statistical model to approximate the value of being in a state resulting from a decision made now.

Policy search has evolved primarily within the computer science community where a policy is always interpreted as some parametric function (\cite{ng2000pegasus}, \cite{peshkin2000learning}, \cite{hu2007evolutionary}, \cite{Robots}, and \cite{mannor2003cross}). Policies might be linear models (where they are referred to as affine policies), locally linear models, and neural networks.  These classes of policies are typically limited to low-dimensional, continuous control problems that arise in engineering.

Below, we identify two major classes of policies that can be used within policy search.  The first includes the parametric functions that are so familiar to the policy search literature.  The second, which is the focus of this research, involves parametrically modified optimization problems which are much more amenable to handling very high-dimensional problems that arise in logistics.

\subsection{Policy Search}
This broad class of algorithmic strategies is based on optimizing a class of parametric functions (policies) by solving the following problem 
	\begin{equation}\label{PolicySearchProb}
		\theta^{*} = \argmin_{\theta \: \in \: \Theta} \mathbb{E}\left[ \sum_{t=0}^T C(S_t,X^\pi_t(S_t|\theta)) \: \big| \: S_0 \right]
	\end{equation}
	where $S_{t+1} = S^M(S_t,x_t, W_{t+1})$, and where $X^{\pi}_t : S_t \times \theta \rightarrow x_t $ is a parametric function whose behavior is determined by a vector of parameters $\theta$.

 A few common examples of parametric policies are the Boltzmann exploration policy, logistic regression, linear decision rules, and neural networks.  Linear decision rules (LDR), also known as affine policies, may be the simplest example of a parametric model. For example, if a decision maker must determine how much water, $x_t$, to release from a reservoir given the capacity of the reservoir, $S_t$, they can use the following decision rule
\begin{equation}
	x_t =\theta_2(S_t - \theta_1) 
\end{equation}
where $\theta_1$ and $\theta_2$ are scalar parameters selected to optimize some predetermined objective function (see \cite{chen2008linear} for more detail). The Boltzmann exploration policy is also parameterized by a positive scalar, $\theta$, where a discrete action, $x$, is determined stochastically with a probability proportional to the estimated value of that action. Formally, the probability of selecting action $x$ given the current state, $S_t$, and parameter $\theta$ is
\begin{equation}
	\mathbb{P}(X_t = x | S_t, \theta) = \frac{ e^{Q(S_t,x) \cdot \theta} }{\sum_{x' \in \mathcal{X}}e^{Q(S_t,x') \cdot \theta}}
\end{equation}
where $Q(S_t, x)$ is the value of decision $x$ given we are in state $S_t$. 

The parameters in each of these previous examples can be selected to solve equation \eqref{PolicySearchProb}. The decision, $x$, typically needs to be very simple - either a discrete set, or a very low dimensional set of continuous actions. This problem is in the same problem class as that posed in classical stochastic search, which is typically written as $ \min_x \mathbb{E} \left[ F(x,W) \right]$. 

In general, there exist two subclasses of algorithms within policy search: gradient-free and gradient-based methods. Policy gradient algorithms approximate the gradient of the expected contribution with respect to the parameter, $\theta$, and then adjust $\theta$ using the approximation to determine optimal values for $\theta$ (see \cite{spall2005introduction}).  These methods build upon the Robbins-Monro algorithm (\cite{robbins1951stochastic}) and are usually some type of adaption of the following steepest descent algorithm:
\begin{equation}\label{eqn:RM}
	\theta^n = \theta^{n-1} + \alpha_{n-1}\nabla_{\theta}\bar{F}(\theta^{n-1},\omega^n)|_{\theta = \theta^{n-1}}
\end{equation}							
where the stepsizes $\alpha_n$ satisfy the following conditions
	\begin{equation}
		\alpha_n > 0, \: \: \: \sum^{\infty}_{n=0} \alpha_n = \infty, \: \: \: \mathbb{E} \left[ \sum^{\infty}_{n=0} (\alpha_n)^2 \right] < \infty, \: \: \text{a.s}.
	\end{equation}	

In general, there are many advantages for using policy gradient methods to optimize parameterized policies. There already exists an extensive literature in stochastic search about stochastic gradient methods (see \cite{spall2005introduction}). More specifically, there exists an extensive and growing literature in reinforcement learning about policy search algorithms and their application to Markov decision processes (see \cite{peters2008reinforcement} and \cite{Robots}). Additionally, many of these algorithms have strong theoretical foundations, convergence guarantees, and are capable of handling multi-dimensional and continuous parameter spaces and situations where some modeling assumptions are unknown. 

Classical policy search is not applicable to high-dimensional, constrained optimization problems that arise in supply chain management, unit commitment problems, or scheduling aircrafts. This problem class has been approached almost exclusively using stochastic programming (a form of lookahead), but this approach produces very large-scale models that can be extremely difficult to solve, and which still suffer from the use of several approximation strategies: limited horizons, replacing multistage with two-stage models, and optimizing over a limited set of scenarios.  Furthermore, as with all policies based on lookahead approximations, stochastic programs, which are basically nonparametric approximations based on sampled representations, ignore obvious structural results which describe the behavior of the policy.

Missing from the literature is the idea of imbedding the parametric approximation within a minimization operator.   We are going to draw on the idea of a parametric cost function approximation, but we are going to combine this with a deterministic lookahead approximation. This is the approach we take in this paper.

\subsection{Policies based on lookahead approximations}		 
Approximate lookahead methods either approximate a \emph{model} of the future, or develop a function that approximates the future value of an action given the current state. We review each of these below. 

\subsubsection{Approximate lookahead models} 
	An optimal policy, $X^{*}_t(\cdot)$, can be defined using the following function
\begin{equation} \label{LA}
	X^*_t(S_t) = \argmin_{x_t} \: \: \left( C(S_t, x_t) + \mathbb{E} \left \{ \min_{\pi \in \Pi} \: \mathbb{E} \left \{ \sum^T_{t' = t+1} C(S_{t'}, X^{\pi}_{t'} (S_{t'})) \big| S_{t+1} \right \} \bigg| S_t, x_t \right \} \right).
\end{equation}
	where $S_{t+1} = S^M(S_t, x_t, W_{t+1})$. We are rarely able to solve equation \eqref{LA} exactly and never for the types of high-dimensional problems we are interested in.  For this reason, we replace the lookahead problem in \eqref{LA} with an approximate lookahead model, which typically draws on five classes of approximating strategies:
\begin{itemize}
	\item {\bf Limiting the horizon} - Here we replace the original horizon $(t,T)$ with a truncated horizon $(t,t+\min\{T,t+H\})$. 
	\item {\bf Stage aggregation} - a stage includes the process of making a decision and then observing exogenous information. Stage aggregation reduces the number of stages in a problem. For example, an $n$-stage problem can be reduced to a 2-stage problem where the modeler makes a decision at time $t$, then observes all remaining exogenous information until time $t + H$ and finally makes the remaining decisions simultaneously.
	\item {\bf Outcome aggregation or sampling} - Here we replace the original outcome space, $\omega \in\Omega$, where $\omega$ is a single realization of $W_1, \ldots, W_T$, with a sampled space $\hat{\Omega}$.
	\item {\bf Discretization} - States, time, and decisions can be discretized in order to make computation tractable. For example, our base model in equation \eqref{objective1} might have time steps of 5 minutes for an energy application, while a lookahead model might use hourly time steps.
	\item {\bf Dimensionality reduction} - In order to simplify the model, a modeler can ignore or remove some variables from the model. These variables may be held fixed as a latent variable in the lookahead models. A common approximation (which we use) is to treat the forecast as fixed in the lookahead model, while in the base model it evolves over time.
\end{itemize}

Lookahead models can take a variety of forms: deterministic lookahead models, also referred to as rolling horizon procedures (\cite{sethi1991theory}) or model predictive control (\cite{camacho2013model}), decision trees (which can be approximated using Monte Carlo tree search) for discrete actions, or stochastic programming models using scenario trees (see \cite{birge2011introduction} and \cite{donohue2006abridged}).

A popular strategy for approximating the exogenous information is to use a sampled set of outcomes that are generally known as scenarios. We call these types of policies \emph{stochastic lookahead models}, but they are also known as scenario optimization and stochastic  programming (see \cite{birge2011introduction}, \cite{dembo1991scenario}, and \cite{mulvey1992stochastic}). To make a distinction between the base model (which is the problem we are trying to solve) and the lookahead model, we use the same notation as in the base model, but we introduce tilde's on all the variables. Each variable carries a double time index $(t,t')$, where $t$ refers to the time at which the lookahead model is formulated in the base model, and $t'$ is the time within the lookahead model. Given this notation, we can write lookahead models as
	\begin{footnotesize}
	\begin{equation}
		X_t^{\text{LA-SP},n}(S_t) = \argmin_{x_t} \: \: \bigg( C(S_t,x_t) + \tilde{\mathbb{E}} \left \{ \min_{(\tilde{x}_{tt'}(\tilde{\omega}),t<t'\leq t+H),\forall \tilde{\omega} \in \tilde{\Omega}^n_t}\tilde{\mathbb{E}}^n \left \{ \sum^{t+H}_{t' = t +1} C(\tilde{S}_{tt'}, \tilde{x}_{tt'}(\tilde{\omega})) \bigg| \tilde{S}_{t,t+1} \right \} \bigg| S_t, x_t \right \} \bigg)
	\end{equation}
	\end{footnotesize}
	where $\tilde{S}_{t, t+1} = S^M(S_t, x_t, \tilde{W}_{t, t+1}(\tilde{\omega}))$ and $\tilde{S}_{t, t'+1} = S^M(\tilde{S}_{tt'}, \tilde{x}_{tt'}, \tilde{W}_{t, t'+1}(\tilde{\omega}))$ where $t' = t+1,..., T-1$. All observations $\tilde{\omega}$ come from the sample set $\tilde{\Omega}^n_t$ and the expectation $\tilde{\mathbb{E}}^n[ \cdot \: | S_t]$ is calculated over the sample set $\tilde{\Omega}^n_t$. We note that all variables indexed by $t$ are $\mathcal{F}_t$-measurable with respect to the information process in the base model.
	These policies are traditionally written as
	\begin{equation}
		X_t^{\text{LA-SP}}(S^n_t) = \argmin_{\tilde{x}_{tt}, (\tilde{x}_{tt'}(\omega), t < t' \leq t + H), \forall \tilde{\omega} \in \tilde{\Omega}_t}  \bigg( c_t x_t + \sum_{\omega \in \tilde{\Omega_t}} \tilde{\mathbb{P}}(\omega) \sum_{t'=t+1}^{t+H} \tilde{c}_{tt'}(\omega) \tilde{x}_{tt'}(\omega) \bigg)
	\end{equation}
\begin{equation*}
	\text{subject to } A_{t}x_{t} \geq b_{t},
\end{equation*}
\begin{equation*}
	A_{tt'}(\tilde{\omega})x_{tt'}(\tilde{\omega})  - B_{t,t'-1}(\tilde{\omega})x_{t,t'-1}(\tilde{\omega})\geq b_{tt'}(\tilde{\omega}) \: \: \tilde{\omega} \in \tilde{\Omega_t}.
\end{equation*}
	These types of policies have received extensive attention in the finance community, particularly the techniques for generating scenarios (see \cite{boender1997hybrid}, \cite{mulvey1992stochastic}, \cite{mulvey1996generating}, \cite{mak1999monte}, \cite{bayraksan2011sequential}, \cite{dupavcova2003scenario}, and \cite{king2012modeling}). These methods are computationally expensive when scaled.They have also received considerable attention for the stochastic unit commitment problem.\cite{ryan2013toward}, \cite{papavasiliou2013multiarea}, \cite{barth2006stochastic}, \cite{zhao2013multi}, \cite{feng2013new}.

Instead of using scenarios a policy may use a point forecast as an approximation for future exogenous information. These are commonly referred to as \emph{deterministic lookahead policies} (see \cite{morari1999model}, \cite{camacho2013model}, and \cite{MP_Mayne}). Using point estimates of future information allows the modeler to use standard deterministic optimization techniques. Below is an example of a deterministic lookahead policy 
\begin{equation}
	X^{\pi}_t (S_t|\theta) = \arg \min_{\tilde{x_{tt}},...,\tilde{x}_{t, t+H}} \sum^{t+H}_{t' = t} \tilde{c}_{tt'} \tilde{x}_{tt'}.
	\end{equation}

Stage aggregation can also be used to approximate the underlying model. The $n$-stage model is often approximated by a two-stage model. The classic two-stage stochastic programming problem, also commonly known as ``recourse" problems (\cite{dantzig1955linear}), is the earliest effort to incorporate uncertainty into linear programs. Formally, the recourse problem is
\begin{equation}
\begin{aligned}
	& \min_{x_0} & & c_0x_0 + \sum_{\omega \in \hat{\Omega}} p(\omega) \min_{x_1(\omega),\omega \in \hat{\Omega}} c_1(\omega)x_1(\omega) \\
\end{aligned}
\end{equation} 
\begin{center}
	\noindent {\bf subject to}
\end{center}
\begin{equation*}
\begin{aligned}
	A_0x_0 & = b_0, \\
	A_1(\omega) x_1(\omega)  - B_0 x_0& = b_0, \\
	x_0 & \leq u_0, \\
	x_1(\omega) & \leq u_1(\omega), \\
	x_0, x_1(\omega) & \geq 0.
\end{aligned}
\end{equation*} 
	
\subsubsection{Value Functions} 
An alternative approach is to use Bellman's equation to approximate the future using value functions (see \cite{putermanmarkov}). A value function, $V_t(S_t)$ is the value of being in a state, $S_t$, and selecting the decision $x_t$. Bellman's equation is formally stated as: 
		\begin{equation}\label{eqn:VA}
	V_{t}(S_{t}) = \min_{x_t} \: \: \left( C(S_t, x_t) + \mathbb{E} \left \{ \min_{\pi \in \Pi} \: \mathbb{E} \left[ \sum^T_{t' = t+1} C(S_{t'}, X^{\pi}_{t'} (S_{t'})) | S_t \right] \bigg| S_t, x_t \right \} \right),
\end{equation}
	where $S_{t+1} = S^M(S_t, X^{\pi}_t, W_{t+1})$.  The minimized expected value of the remaining cumulative contributions is equivalent to the expected value of $V_{t+1}(S_{t+1})$. Therefore, equation \eqref{eqn:VA} can be rewritten as
	\begin{equation} 
			V_{t}(S_{t}) = \min_{x_t} \: \: \left( C(S_t, x_t) + \mathbb{E} \bigg [ V_{t+1}(S_{t+1}) \bigg | S_t \bigg ] \right)
		\end{equation}
		where $S_{t+1} = S^M (S_t, x_t, W_{t+1})$.
	Computationally, expectations can be very difficult, especially for high dimensional random variables. Computational burdens regarding expectation can be circumvented by approximating the value function based on the post-decision state, $S^x_t$ which is the state of the system immediately after the decision maker has made the decision, $x_t$. This removes the need for the imbedded expectation. If the modeler is approximating around the post-decision then the policy takes the form
\begin{equation}
	X^{\pi}_t(S_t | \theta) = \arg \min_{x \in \chi_t} \left(C(S_t, x) + \bar{V}^x_{t}(S^x_{t}|\theta) \right).
\end{equation}
	
The modeler must determine how to approximate the value function. There is by now an extensive literature that addresses approximating value functions known as approximate (or adaptive) dynamic programming, or reinforcement learning (\cite{ADP_Powell}, \cite{bertsekas2011dynamic}, and \cite{sutton1998reinforcement}).  Methods such as Stochastic Dual Dynamic Programming (SSDP), linear approximations, and piecewise linear approximation scale to higher dimensions (see \cite{simao2009approximate}, \cite{godfrey2002adaptive}, \cite{topaloglu2006dynamic}, \cite{pereira1991multi}, \cite{shapiro2014lectures}, \cite{shapiro2011analysis}, and \cite{philpott2008convergence}).  In this setting, value functions can be approximated using Benders cuts. This strategy has been widely used for solving multistage stochastic programming problems using Benders cuts (\cite{girardeau2014convergence}, \cite{sen2014multistage}, \cite{pereira1991multi}, \cite{shapiro2009lectures}, \citet{birge2011introduction}) or piecewise linear, separable value function approximations (\cite{topaloglu2006dynamic}, \cite{simao2009approximate}, \cite{ADP_Powell}). This literature has produced a wide range of approximations that produce effective (if not optimal) policies. 

As of this writing, forecasts (which evolve over time) have never been modeled explicitly as part of the state variable (as they should).  When forecasts are left out of the state variable, that means they are being treated as latent variables, which means that the dynamic program is actually a lookahead model which would have to be solved again when the forecasts change.  By contrast, it is straightforward to incorporate forecasts within lookahead models, since these are, in fact, reoptimized as we step forward in time (and update the forecasts). 

Other work tuned for complex dynamic resource allocation problems that arise in transportation and logistics has been developed using the language of approximate dynamic programming (\cite{topaloglu2006dynamic}, \cite{simao2009approximate}, \cite{ADP_Powell}). We note only that any method based on approximating a value function, whether with Benders cuts or statistical machine learning, is a form of approximate dynamic programming.	 

\section{The Parametric Cost Function Approximation}
We extend the concept of policy search to include parameterized optimization problems. The parametric Cost Function Approximation (CFA) imitates the structural simplicity of deterministic lookahead models and myopic policies, but allows more flexibility by adding tunable parameters. 

\subsection{Basic Idea}
Since the idea of a parametric cost function approximation is new, we begin by outlining the general strategy, and then demonstrate how to apply it for our energy storage problem. We propose using parameterized optimization problems such as 
	\begin{equation}
		X^{\pi}_t(S_t|\theta) = \argmin_{x_t \in \mathcal{X}^\pi(\theta)} \: \bigg\{C(S_t,x_t) + \sum_{f \in \mathcal{F}} \theta^c_f \phi_f (S_t, x_t)  : A_tx_t = \bar{b}^{\pi}_t(\theta^b) \bigg\}
	\end{equation}
as a type of parameterized policy. Here the index $\pi$ signifies the structure of the modified set of constraints, $\theta^c$ is the vector of cost function parameters, $\theta^b$ is the vector of constraint parameters, and $\phi_f$ are the basis functions corresponding to features $f \in \mathcal{F}$. 

Parametric terms can also be added to the cost function or constraints of a myopic or deterministic lookahead model. In the following example, parameters have been added as an error correction term to the contribution function as well as to the model constraints.
\begin{equation}
	X^{\pi}_t(S_t | \theta) = \argmin_{x_t \in \mathcal{X}_t} \: \: \left( C(S_t, x_t) + \sum_f \theta^c_f \phi_f (S_t, x_t) \right)
\end{equation}
$$ \text{\bf subject to } A_tx_t = b_t + D\theta^b$$ where $D$ is a scaling matrix. We emphasize that the cost correction term should not be confused as a value function approximation, because we make no attempt to approximate the downstream value of being in a state.

\subsection{A hybrid Lookahead-CFA policy}
There are many problems that naturally lend themselves to a lookahead policy (for example, to incorporate a forecast or to produce a plan over time), but where there is interest in making the policy more robust than a pure deterministic lookahead using point forecasts.  For this important class (which is the problem we face), we can create a hybrid policy where a deterministic lookahead has parametric modifications that have to be tuned using policy search. When parameters are applied to the constraints it is possible to incorporate easily recognizable problem structure.  For example, a supply chain management problem can handle uncertainty through buffer stocks, while an airline scheduling model might handle stochastic delays using schedule slack.  A grid operator planning energy generation in the future might schedule reserve capacity to account for uncertainty in forecasts of demand, as well as energy from wind and solar.  As with all policy search procedures, there is no guarantee that the resulting policy will be optimal unless the parameterized space of policies includes the optimal policy. However, we can find the optimal policy within the parameterized class, which may reflect operational limitations.  We note that while parametric cost function approximations are widely used in industry, optimizing within the parametric class is not done.

\subsection{Parameterization Structure}
	Parametric terms can be appended to existing  constraints, new parameterized constraints can be added to the existing model, or a combination of the two can be done. Often the problem setting will influence how policy constraints should be parameterized. Consider the energy storage problem where a manager must satisfy the power demand of a small building. The manager has a stochastic supply of renewable energy at no cost, an unlimited supply from the main power grid at a stochastic price, and access to a local rechargeable storage device. Every period the manager must determine what combination of energy sources to use to satisfy the power demand, how much energy to store, and how much to sell back to the grid. Given the manager has access to point forecasts of future exogenous information he or she can use the following lookahead policy to determine how to allocate their energy.  
	
\begin{equation} \label{energy_LA}
	X^{\text{D-LA}}_t(S_t) = \argmin_{x_{tt}} \: \: \left( C(S_{t}, x_{t}) + \: \left[ \sum^T_{t' = t+1} \tilde{c}_{tt'} \tilde{x}_{tt'} \right] \right)
\end{equation}
where $\bar{S}_{t,t'+1} = \tilde{S}^M(\bar{S}_{tt'}, \tilde{x}_{tt'}, \bar{W}_{t,t'+1})$. It is important to note that if the contribution function, transition function, and constraints of $X^{\pi}(\cdot)$ are linear, this policy can be expressed as the following linear program. 
\begin{equation}\label{PCFA}
		X^{\text{D-LA}}_t(S_t) = \argmax_{x_{t}} \: \{c_t^Tx_t : A_t x_t \leq b_t, \:\: x_t \geq 0  \}.
\end{equation}
There are different ways to parameterize the previous policy, but since all the uncertainty in our problem is restricted to the right hand side constraints, we will only parameterize the vector $b_t$. Once parameterized our policy becomes
\begin{equation}
	X^{\text{LA-CFA}}_t(S_t|\theta) = \argmax_{x_{t}} \: \{c_t^Tx_t : A_t x_t \leq b_t(\theta)  \}
\end{equation}
	where $\theta$ is a vector of tunable parameters. Parametric modifications can be designed specifically to capture the structure of a particular policy. The idea to use buffers and inventory constraints to manage storage is intuitive and easily incorporated into a deterministic lookahead. In the previous energy storage problem a lower buffer guarantees the decision maker will always have access to some stored energy. Conversely, an upper threshold will make sure some storage space remains in the battery in order to capitalize on unexpected gusts of wind (for example). For the energy storage problem, we represent the amount of energy stored in the battery as the variable $R_t$ and the approximated future amount of energy in storage at time $t'$ given the information available at time $t$ as $R_{tt'}$.  Thus, 
\begin{equation}
	\theta^L \leq R_{tt'} \leq \theta^U \text{ for } t' > t.
\end{equation}
Although it can greatly increase the parameter space, the upper and lower bounds can can also depend on $(t'-t)$
\begin{equation}\label{TimeIN}
	\theta^L_{t'-t} \leq R_{tt'} \leq \theta^U_{t'-t} \text{ for } t' > t.
\end{equation}
The resulting modified deterministic problem is no harder to solve than the original deterministic problem (where $\theta^L = 0$ and $\theta^U = R^{max}$).  We now have to use policy search techniques to optimize $\theta$. Below we suggest different ways of parameterizing the right hand side adjustment.

\subsubsection{Lookup table in time}
Policy parameterizations come in a variety of forms. A simple form is a lookup table indexed by time as in equation (\ref{TimeIN}). Though it may be simple, a lookup table model for $\theta$ means that the dimensionality increases with the horizon which can complicate the policy search process.  

In the energy storage example, $f^E_{tt'}$ represents the forecast of the amount of renewable energy available at time $t'$ given the information available at time $t$. A policy maker may use the parametrization $\theta_{t'-t} \cdot f^E_{tt'}$ to intentionally overestimate or underestimate the amount of future renewable energy. The policy maker may set $\theta_{t'-t} \leq 1$ to make the policy more robust and avoid the risk of running out of energy. 

This type of parameterization is not limited to just modifying the point forecast of exogenous information. If the modeler has sufficient information such as the cumulative distribution function of $E_{t'}$, $F^{E_{t'}}(\cdot)$, he or she may even exchange the point forecast $f^E_{tt'}$ with the quantile function
	\begin{equation}
		Q^E_{t'}(\theta_{t'-t}) = \inf \left\{w\in {\mathbb  {R}}: \theta_{t'-t} \leq F^{E_{tt'}}(w)\right\}.
	\end{equation}
In this case $\theta_{t'-t'}$ is still a parameter of the policy and determines how aggressively or passively the policy stores energy. The lookup table in time parameterization is best if the relationship between parameters in different periods is unknown.

\subsubsection{Parametric model}
Instead of having an adjustment $\theta_\tau = \theta_{t'-t}$ for each time $t+\tau$ in the future, we can use instead a parametric function of $\tau$, which reduces the number of parameters that we have to estimate. For example, we might use the parametric adjustment:
\begin{equation}
	\theta^L\cdot e^{\alpha \cdot (\tau)} \leq R_{tt'} \leq \theta^U \cdot e^{\beta \cdot (\tau)} \text{ for } t' > t \text{ and } \alpha, \beta \in \mathbb{R}.
\end{equation}

These parametric functions of time can also be used to directly modify the lookahead model's forecasts. For example, in the energy storage example, the policy maker may use the parameterization $f^E_{tt'} \cdot \theta_1 e^{\theta_2\cdot (t'-t)}$ to replace the forecasted amount of future renewable energy, $f^E_{tt'}$.

\section{Determining parameters for the CFA}
Policy Function Approximations (PFAs) and Cost Function Approximations (CFAs) are structurally different in the sense one uses analytic functions and the other uses parameterized optimization problems, respectively, to make decisions. However, they are both subclasses of the same general class of parameterized policies, $X^{\pi}_t : \mathcal{S}_t \times \Theta \rightarrow \mathcal{X}_t$, and their optimal parameterization, $\theta^*$, can be found by solving
\begin{equation}\label{eqn:param_prob}
	\theta^* = \argmax_{\theta \in \Theta} \: F(\theta),
\end{equation} 
where
\begin{equation}\label{cum_reward}
	F(\theta) = \mathbb{E}\left[ \sum_{t=0}^T C(S_t,X^\pi_t(S_t|\theta)) \: \big| \: S_0 \right]
\end{equation}
and $S_{t+1} = S^M(S_{t}, X^\pi_t(S_t|\theta), W_{t+1})$. If $F(\cdot)$ is well defined, finite valued, convex, and continuous at every $\theta$ in the nonempty, closed, bounded, and convex set $\Theta \subset \mathbb{R}^n$, then an optimal $\theta^* \in \Theta$ exists. It is possible to use the iterative algorithm described in equation \eqref{eqn:RM} to find $\theta^*$, (see \cite{robbins1951stochastic}). If $W_{t}$ is a stochastic process adapted to the filtration $(\mathcal{F}_t)_{t\geq0}$, there exists a stochastic subgradient, $g^n \in \partial_{\theta}F(\theta^{n-1})$, and the following assumptions are satisfied.
	\begin{enumerate}
		\item $\mathbb{E}\bigg[g^{n+1} \cdot (\bar{\theta}^{n} - \theta^*) \bigg| \mathcal{F}^n\bigg] \geq 0$.
		\item $|g^n| \leq B_g$.
		\item For any $\theta$ where $|\theta - \theta^*| > \delta, \delta > 0$, there exists $\epsilon > 0$ such that $\mathbb{E}[g^{n+1}|\mathcal{F}^n] > \epsilon$.
	\end{enumerate}
This is true regardless if the policy $X^{\pi}_t(\cdot|\theta)$ is a CFA or PFA. There are several ways to generate stochastic subgradients that satisfy the previous conditions. If the cumulative reward of a single sample path, $\bar{F}(\cdot, \omega)$, is convex and differentiable for every $\omega \in \Omega$ and $\theta$ is an interior point of $\Theta$ , then the gradient, $\nabla_{\theta}\bar{F}$, of
			\begin{equation}\label{sample_cum_reward}
				\bar{F}(\theta, \omega) = \sum_{t=0}^T C\bigg(S_{t}(\omega),X^\pi_{t}(S_{t}(\omega)|\theta) \bigg|\theta \bigg)
			\end{equation}
where $S_{t+1}(\omega) = S^M(S_t(\omega), X^{\pi}_t(S_t(\omega)), W_{t+1}(\omega))$ can be used as an appropriate stochastic gradient (see \cite{strassen1965existence}). This subgradient can be calculated recursively and is described in the following proposition.		
	
\newtheorem{Title2}{Proposition}
\begin{Title2}
	Assume $\bar{F}(\cdot, \omega)$ is convex for every $\omega \in \Omega$, $\theta$ is an interior point of $\Theta$, and $F(\cdot)$ is finite valued in the neighborhood of $\theta$, then 
	\begin{equation*}
		\nabla_{\theta}F(\theta) = \mathbb{E}[\nabla_{\theta}\bar{F}(\theta, \omega)]
	\end{equation*}
	where 
	\begin{equation*}
		\nabla_{\theta}\bar{F} = \bigg(\frac{\partial C_0}{\partial X_0} \cdot \frac{\partial X_0}{\partial \theta} \bigg) + \sum^T_{t'=1} \bigg[ \bigg( \frac{\partial C_{t'}}{\partial S_{t'}} \cdot \frac{\partial S_{t'}}{\partial \theta} \bigg) + \bigg( \frac{\partial C_{t'}}{\partial X_{t'}(S_t | \theta)} \cdot \bigg( \frac{\partial X_{t'}(S_t | \theta)}{\partial S_{t'}} \cdot \frac{\partial S_{t'}}{\partial \theta} + \frac{\partial X_{t'}(S_t | \theta)}{\partial \theta} \bigg) \bigg) \bigg], \\
        \end{equation*}
        and
        \begin{equation*}
            	\frac{\partial S_{t'}}{\partial \theta} = \frac{\partial S_{t'}}{\partial S_{t'-1}} \cdot \frac{\partial S_{t'-1}}{\partial \theta} + \frac{\partial S_{t'}}{\partial X_{t'-1}(S_{t-1}|\theta)} \cdot \bigg[ \frac{\partial X_{t;-1}(S_{t-1}|\theta)}{\partial S_{t'-1}} \cdot \frac{\partial S_{t'-1}}{\partial \theta} + \frac{\partial X_{t'-1}(S_{t-1}|\theta)}{\partial \theta} \bigg].
        \end{equation*}   
\end{Title2}

\begin{proof}
If $\bar{F}(\cdot, \omega)$ is convex for every $\omega \in \Omega$, $\theta$ is an interior point of $\Theta$, and $F(\cdot)$ is finite valued in the neighborhood of $\theta$, then by theorem 7.47 of \cite{SP_Shapiro} 
	\begin{equation*}
		\nabla_{\theta} \mathbb{E} F(\theta,W) = \mathbb{E} \nabla_{\theta} F(\theta,W).
	\end{equation*}
	Applying the chain rule, we find:
	\begin{equation*}
	\begin{split}
	 \nabla_{\theta}\bar{F} &= \frac{\partial}{\partial \theta} \bigg[ C_0(S_0, X_0(S_0|\theta )) + \sum^T_{t'=1} C(S_{t'},X_{t'}(S_{t'}|\theta )) \bigg] \\
                        		&= \frac{\partial}{\partial \theta} C_0(S_0, X_0(S_0|\theta)) + \frac{\partial}{\partial \theta} \bigg[  \sum^T_{t'=1} C(S_{t'},X_{t'}(S_{t'}|\theta)) \bigg] \\
                        		&= \bigg(\frac{\partial C_0}{\partial X_0} \cdot \frac{\partial X_0}{\partial \theta} \bigg) + \bigg[  \sum^T_{t'=1} \frac{\partial}{\partial \theta} \: C(S_{t'},X_{t'}(S_{t'}|\theta)) \bigg] \\
                        		&= \bigg( \frac{\partial C_0}{\partial X_0} \cdot \frac{\partial X_0}{\partial \theta} \bigg) + \sum^T_{t'=1} \bigg[ \bigg( \frac{\partial C_{t'}}{\partial S_{t'}} \cdot \frac{\partial S_{t'}}{\partial \theta} \bigg) + \bigg( \frac{\partial C_{t'}}{\partial X_{t'}(S_{t'}|\theta)} \cdot \frac{\partial X_{t'}(S_{t'}|\theta)}{\partial \theta} \bigg) \bigg] \\
                        		&= \bigg(\frac{\partial C_0}{\partial X_0} \cdot \frac{\partial X_0}{\partial \theta} \bigg) + \sum^T_{t'=1} \bigg[ \bigg( \frac{\partial C_{t'}}{\partial S_{t'}} \cdot \frac{\partial S_{t'}}{\partial \theta} \bigg) + \bigg( \frac{\partial C_{t'}}{\partial X_{t'}(S_{t'}|\theta)} \cdot \bigg( \frac{\partial X_{t'}(S_{t'}|\theta)}{\partial S_{t'}} \cdot \frac{\partial S_{t'}}{\partial \theta} + \frac{\partial X_{t'}(S_{t'}|\theta)}{\partial \theta} \bigg) \bigg) \bigg], \\
                        \end{split}
                        \end{equation*}
                        where
                         \begin{equation*}
            	\frac{\partial S_{t'}}{\partial \theta} = \frac{\partial S_{t'}}{\partial S_{t'-1}} \cdot \frac{\partial S_{t'-1}}{\partial \theta} + \frac{\partial S_{t'}}{\partial X_{t'-1}(S_{t-1}|\theta)} \cdot \bigg[ \frac{\partial X_{t;-1}(S_{t-1}|\theta)}{\partial S_{t'-1}} \cdot \frac{\partial S_{t'-1}}{\partial \theta} + \frac{\partial X_{t'-1}(S_{t-1}|\theta)}{\partial \theta} \bigg].
        \end{equation*} 
\end{proof}
		
\subsection{Computing the gradient of linear cumulative reward}
If the objective function in equation (\ref{cum_reward}) is a linear function of the decisions, $x_t$, the parametric CFA policy, $X^{\pi}_t(S_t|\theta)$, which determines the decision, $x_t$, can be written as the following linear program:
    \begin{equation*}
    	X^{\pi}_t(S_t|\theta) = \argmax_{\tilde{x}_{tt}} \: \: \sum_{t' = t}^T c_{tt'} \tilde{x}_{tt'}
    \end{equation*}
where $A\tilde{x_t} \leq b(\tilde{W}_{tt'},\theta)$ and $\tilde{x}_t^T = [\tilde{x}_{t,t},...,\tilde{x}_{t,T}]$.  Point estimates, $\tilde{W}_{tt'}$, can be used to approximate exogenous information. If this policy is written as a linear program where the state and approximated exogenous information is only in the right hand side constraints, $b(\tilde{W}_{tt'}, \theta)$, then a subgradient can be calculated recursively and is described in the following proposition:

\newtheorem{Title4}{Proposition}
\begin{Title4}
Given $\bar{F}(\cdot, \omega)$ is convex for every $\omega \in \Omega$, $\theta$ is an interior point of $\Theta$, and the contribution function $C(x)$ is a linear function of $x$, the transition function $S_t = S^M(S_{t-1},x_{t-1},W_{t})$ is linear in $S_{t-1}$ and $x_{t-1}$, $F(\cdot)$ is finite valued in the neighborhood of $\theta$, and the policy, $X^{\pi}_t(S_t|\theta)$ is defined as
\begin{equation}
	X^{\pi}_t(S_t|\theta) = \argmax_{\tilde{x}_{t,t}} \: \: \sum_{t' = t}^T c_{tt'}^T \tilde{x}_{tt'}
	\label{eqn:LP_main}
\end{equation}
where $A\tilde{x_t} \leq b(\tilde{W}_{tt'},\theta)$,  $B_t$ is the basis matrix corresponding to the basic variables for the optimal solution of (\ref{eqn:LP_main}), and $\tilde{x}_t^T = [\tilde{x}_{t,t},...,\tilde{x}_{t,T}]$. Then 
	\begin{equation*}
		\nabla_{\theta}F(\theta) = \mathbb{E}[\nabla_{\theta}\bar{F}(\theta, \omega)]
	\end{equation*}
	where 
	\begin{equation}
		\nabla_{\theta} \bar{F}(\theta, \omega) = \sum^T_{t=1} \bigg( \nabla_{\theta} b(\theta, S_t) + \nabla_{S_t} b(\theta, S_t) \cdot \nabla_{\theta} S_t \bigg)^T \cdot \bigg( B^{-1}_t \bigg)^T \cdot c_t, \\
	\end{equation}
	\begin{equation}
		\nabla_{\theta}S_t = \nabla_{S_{t-1}} S^M(S_{t-1},x_{t-1},W_t) \cdot \nabla_{\theta}S_{t-1} + \nabla_{x_{t-1}} S^M(S_{t-1},x_{t-1},W_t) \cdot \nabla_{\theta} x_{t-1}.
	\end{equation}

\end{Title4}

\begin{proof}
If $\bar{F}(\cdot, \omega)$ is convex for every $\omega \in \Omega$, $\theta$ is an interior point of $\Theta$, and $F(\cdot)$ is finite valued in the neighborhood of $\theta$ then 
	\begin{equation*}
		\nabla_{\theta}F(\theta) = \mathbb{E}[\nabla_{\theta}\bar{F}(\theta, \omega)].
	\end{equation*}
	If the contribution function $C(x)$ is a linear function of $x$, the transition function $S_t = S^M(S_{t-1},x_{t-1},W_{t})$ is linear, and the policy, $X^{\pi}_t(S_t|\theta)$ is defined as
\begin{equation*}
	X^{\pi}_t(S_t|\theta) = \argmax_{\tilde{x}_{t,t}} \: \: \sum_{t' = t}^T c_{tt'}^T \tilde{x}_{tt'}
\end{equation*}
where $A\tilde{x_t} \leq b(S_t,\theta)$ and $\tilde{x}_t^T = [\tilde{x}_{t,t},...,\tilde{x}_{t,T}]$. Then
	\begin{equation}
	\begin{split}
		\nabla_{\theta} \bar{F}(\theta, \omega)	& = \nabla_{\theta} \bigg[ \sum^T_{t=1} c^T_t x_t(S_t|\theta) \bigg] \\
			&= \sum^T_{t=1} \bigg[ \nabla_{\theta} \bigg( c^T_t x_t(S_t|\theta) \bigg) \bigg] \\
			&= \sum^T_{t=1} \bigg[ \nabla_{\theta} x_t(S_t|\theta)^T \bigg] \cdot c_t \\
			&= \sum^T_{t=1} \nabla_{\theta} \bigg[ B^{-1}_t \cdot b(S_t|\theta) \bigg]^T \cdot c_t \\
			&= \sum^T_{t=1} \nabla_{\theta} \bigg[ b(\theta, S_t) \cdot (B^{-1}_t)^T \bigg] \cdot c_t \\
			&= \sum^T_{t=1} \bigg( \nabla_{\theta} b(\theta, S_t)^T \bigg) \cdot \bigg( B^{-1}_t \bigg)^T \cdot c_t \\
			&= \sum^T_{t=1} \bigg( \partial_{\theta} b(\theta, S_t) + \partial_{S_t} b(\theta, S_t) \cdot \partial_{\theta} S_t \bigg)^T \cdot \bigg( B^{-1}_t \bigg)^T \cdot c_t. \\
	\end{split}
	\end{equation} 
\end{proof}

\subsection{The CFA Gradient Algorithm}
The ability to calculate an unbiased estimator of $\nabla_{\theta}F(\theta)$ allows us to use stochastic approximation techniques to determine the optimal parameters, $\theta$, of the CFA policy, $X^{\pi}_t(\cdot | \theta)$. Below is the iterative algorithm we use to tune our CFA policies. 
   \begin{algorithm}
   \caption{CFA Gradient Algorithm}
   \label{CFA_Grad}
   \begin{algorithmic}[1]
    	\State Initialize $\theta^0$, $N$, and $k$:
					\For{ n = 1, 2, 3, ..., N}
						\State Generate a trajectory $\omega^{n}$ where 
							$$S^n_{t+1}(\omega^n) = S^M(S^n_t(\omega^n), X_t^\pi(S^n_t(\omega^n)|\theta^{n-1}),W_{t+1}(\omega^n))$$
						\State Compute the gradient estimator, $\nabla_{\theta}\bar{F}(\theta^{n-1},\omega^n)$							
							
						\State Update policy parameters, $\theta$
							\begin{equation}
								\theta^n = \theta^{n-1} + \alpha_{n-1}\nabla_{\theta}\bar{F}(\theta^{n-1},\omega^n)|_{\theta = \theta^{n-1}}
							\end{equation}
					\EndFor
    \end{algorithmic}
    \end{algorithm}
where the stepsizes $\alpha_n$ satisfy the following conditions
	\begin{equation*}
		\alpha_n > 0, \: \: \: \sum^{\infty}_{n=0} \alpha_n = \infty, \: \: \: \mathbb{E} \left[ \sum^{\infty}_{n=0} (\alpha_n)^2 \right] < \infty, \: \: \text{a.s}.
	\end{equation*}	
If $F(\cdot)$ is continuous and finite valued  in the neighborhood of every $\theta$, in the nonempty, closed, bounded, and convex set $\Theta \subset \mathbb{R}^n$ such that $\bar{F}(\cdot,\omega)$ is convex for every $\omega \in \Omega$ where $\theta$ is an interior point of $\Theta$, then
\begin{equation*}
	\lim_{n \rightarrow \infty} \: \theta^n \longrightarrow \theta^* \: \: \text{a.s.}
\end{equation*}
Although any step size rule that satisfies the previous conditions will guarantee asymptotic convergence, we prefer parameterized rules that can be tuned for quicker convergence rates. In practice the number of updating iterations, $N$, is a finite number. Therefore, we limit our evaluation of the algorithm to how well it does within $N$ iterations. Algorithm \ref{CFA_Grad} can be described as a policy, $\theta^{\pi}(S^n)$, with a state variable, $S^n = \theta^n$ plus any parameters needed to compute the stepsize policy, and where $\pi$ describes the structure of the stepsize rule. If $\theta^{\pi,n}$ is the estimate of $\theta$ using stepsize policy $\pi$ after $n$ iterations, then our goal is to find the policy that produces the best performance (in expectation) after we have exhausted our budget of $N$ iterations.  Thus, we wish to solve
	\begin{equation}
		\max_{\pi} \mathbb{E} \bar{F}(\theta^{\pi,N},W).
	\end{equation}
Our goal is now to find the best stepsize rule that maximizes terminal value within $N$ iterations. For our numerical example we use the adaptive gradient algorithm, ADAGRAD, as our step size rule ( \cite{duchi2011adaptive}). ADAGRAD modifies the individual step size for the updated parameter, $\theta$, based on previously observed gradients. The step size $\alpha_n$  as 
	\begin{equation}
		\alpha_n =  \frac{\eta}{\sqrt{G_{t} + \epsilon}}
	\end{equation}
	where $\eta$ is a scalar learning rate, $G \in \mathbb{R}^{d \times d}$ is a diagonal matrix where each diagonal element is the sum of the squares of the gradients with respect to $\theta$ up to the current iteration $n$, while $\epsilon$ is a smoothing term that avoids division by zero. For our simulations we set $\eta = .1$.
This method is applied to a numerical example in the following section.

\section{An Energy Storage Application}
\begin{figure}[h]
    \centering
    \includegraphics[width=0.6\textwidth]{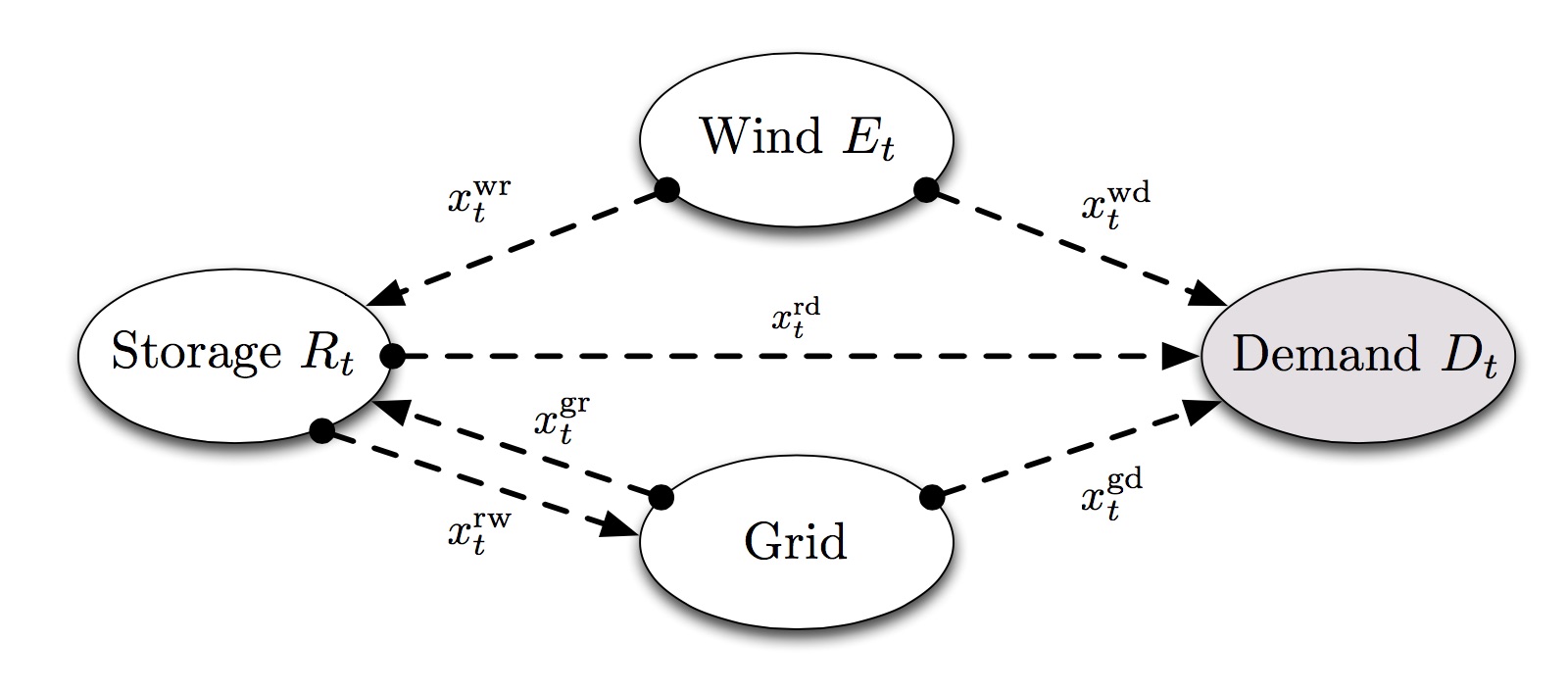}
    \caption{Energy system schematics}
    \label{fig:system}
\end{figure}

To illustrate the capability of the parametric cost function approximation, we use it to solve a time-dependent stochastic energy storage problem where we have access to rolling forecasts of varying quality.  We show how we can use parametrically modified deterministic lookahead models to produce robust policies that work better under uncertainty than a standard deterministic lookahead. In our setting a smart grid manager must satisfy a recurring power demand with a stochastic supply of renewable energy, limited supply of energy from the main power grid at a stochastic price, and access to a local rechargeable storage devices. This system is graphically represented in Figure \ref{fig:system}. 

Every hour the manager must determine what combination of energy sources to use to satisfy the power demand, how much energy to store, and how much to sell back to the grid. The state variable at time $t$, $S_t$, includes the level of energy in storage, $R_t$, the amount of energy available from wind, $E_t$, the spot price of electricity, $P_t$, the demand $D_t$, and the energy available from the grid $G_t$ at time $t$. The state of the system can be represented by the following five dimensional vector,
	\begin{equation}
		S_t = (R_t, E_t, P_t, D_t,G_t)
	\end{equation}
where $R_t \in [0, R_{\max}]$ is the level of energy in storage at time $t$. The demand ($D_t$) has a deterministic seasonal structure
\begin{equation}
	D_t = \lfloor \max \{0 ,100 - 50 \sin \bigg( \frac{5 \pi t}{T} \bigg) \} \rfloor.
\end{equation}

At the beginning of every period $t$ the manager must combine energy from the following sources to satisfy the demand, $D_t$:
\begin{enumerate}
	\item Energy currently in storage (represented by a decision $x_t^{rd}$);
	\item Newly available wind energy (represented by a decision $x_t^{wd}$);
	\item And energy from the grid (represented by a decision $x_t^{gd}$).
\end{enumerate}
Additionally, the manager must decide how much renewable energy to store, $x_t^{wr}$, how much energy to sell to the grid at price $P_t$, $x_t^{rg}$, and how much energy to buy from the grid and store, $x_t^{gr}$. The manager's decision is defined as the following vector
	\begin{equation}\label{eqn:dec}
		x_t = (x^{wd}_t, x^{gd}_t, x^{rd}, x^{wr}_t, x^{gr}_t, x^{rg}_t)^T \geq 0
	\end{equation}
given the following constraints:
	\begin{equation}\label{eqn:con}
		\begin{matrix}
		x^{wd}_t & + & \beta^dx^{rd}_t & + & x^{gd}_t & \leq 		& D_t, \\
		~ & ~ & x^{gd}_t & + & x^{gr}_t & \leq 	& G_t, \\
		~ & ~ & x^{rd}_t & + & x^{rg}_t & \leq 	& R_t, \\
		~ & ~ & x^{wr}_t & + & x^{gr}_t & \leq 	& R_{\max} - R_t, \\
		~ & ~ & x^{wr}_t & + & x^{wd}_t & \leq 	& E_t, \\
		~ & ~ & x^{wr}_t & + & x^{gr}_t & \leq 	& \gamma^c, \\
		~ & ~ & x^{rd}_t & + & x^{rg}_t & \leq 	& \gamma^d \\
		\end{matrix}
\end{equation}
where $\gamma^c$ and $\gamma^d$ are the maximum amount of energy that can be charged or discharged from the storage device. Typically, $\gamma^c$ and $\gamma^d$ are the same. 

The transition function, $S^M(\cdot)$,  explicitly describes the relationship between the state of the model at time $t$ and $t+1$,
\begin{equation*}
	S_{t+1} = S^M(S_t, x_t, W_{t+1})
\end{equation*}
where $W_{t+1} = (E_{t+1}, P_{t+1}, D_{t+1})$ is the exogenous information revealed at $t+1$. In the problem only storage is carried from one period to the next. Price, $P_t$, and demand, $D_t$, do not depend on the past. The relationship of storage levels between periods is defined as
	\begin{equation}\label{eqn:trans}
		R_{t+1} = R_{t} - x^{rd}_t + \beta^c x^{wr}_t + \beta^c x^{gr}_t - x^{rg}_t
	\end{equation}
where $\beta^c \in (0, 1)$ and $\beta^d \in (0, 1)$, are the charge and discharge efficiencies. For a given state $S_t$ and decision $x_t$, we define:
\begin{equation}
	C(S_t, x_t) = P_t \cdot (x_t^{wd} + \beta^dx^{rd} + x^{gd} + \beta^d x^{rg}_t - x^{gr}_t - x^{gd}_t) - C^{\text{penalty}} \cdot \left( D_t - x_t^{wd} - \beta^dx^{rd} - x^{gd} \right)
\end{equation}
where $C^{\text{penalty}}$ is the penalty of not satisfying demand and $C(S_t,x_t)$ is the profit realized at $t$ given the current state is $S_t$ and the decision is $x_t$. The objective is to find the policy $\pi$ that solves
\begin{equation} \label{objective}
\begin{aligned}
	&\max_{\pi \in \Pi} \: \: \mathbb{E}^{\pi} \bigg[ \sum^T_{t=0} C(S_t, X^{\pi}_t(S_t)) \: \bigg| \: S_0  \bigg] \\
	&\text{subject to \eqref{eqn:dec} - \eqref{eqn:trans} for $t \in [1,T]$}.
\end{aligned}
\end{equation}

\subsection{Energy generation model}
\begin{figure}[h]
    \centering
    \includegraphics[width=0.6\textwidth]{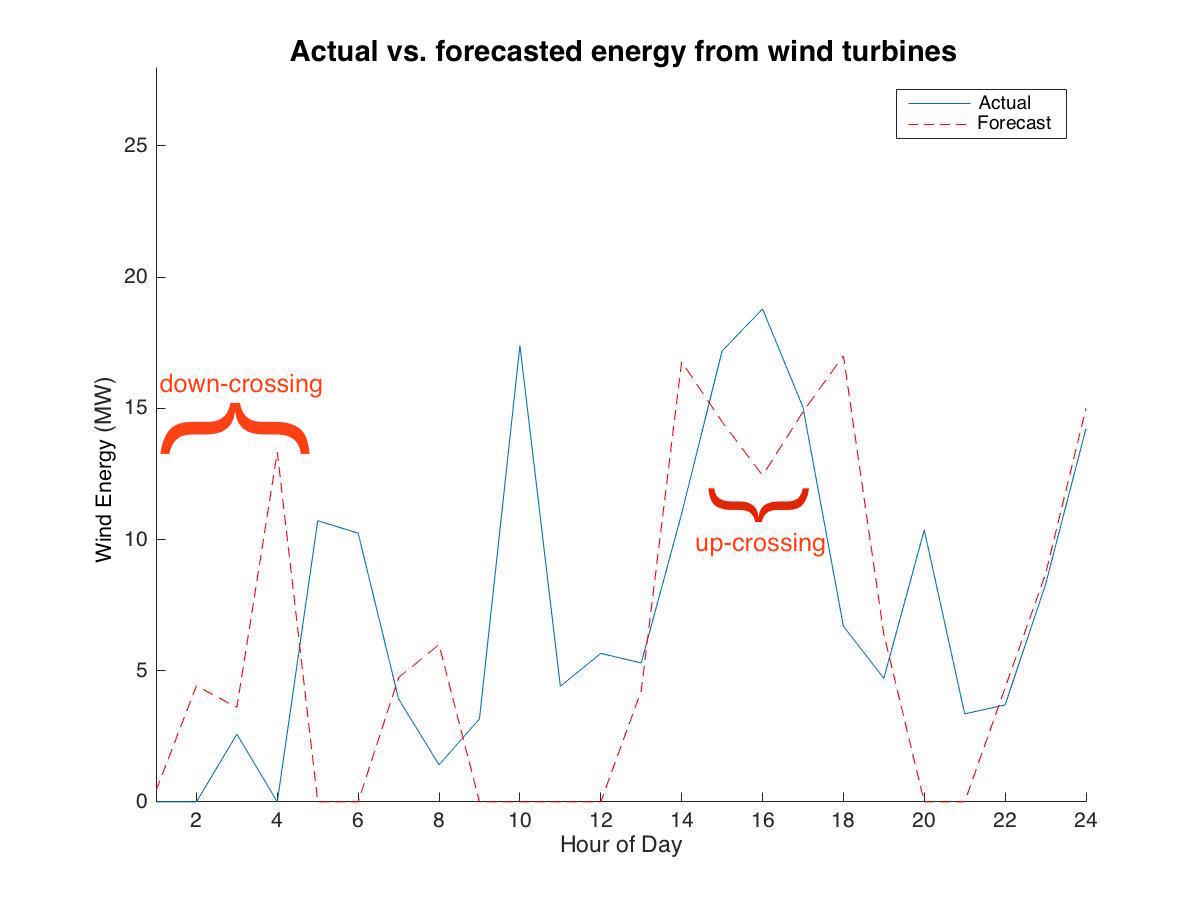}
    \caption{Sample path of renewable energy ($E_t$)}
    \label{fig:renewables}
\end{figure}

\begin{figure}[h]
    \centering
    \includegraphics[width=0.6\textwidth]{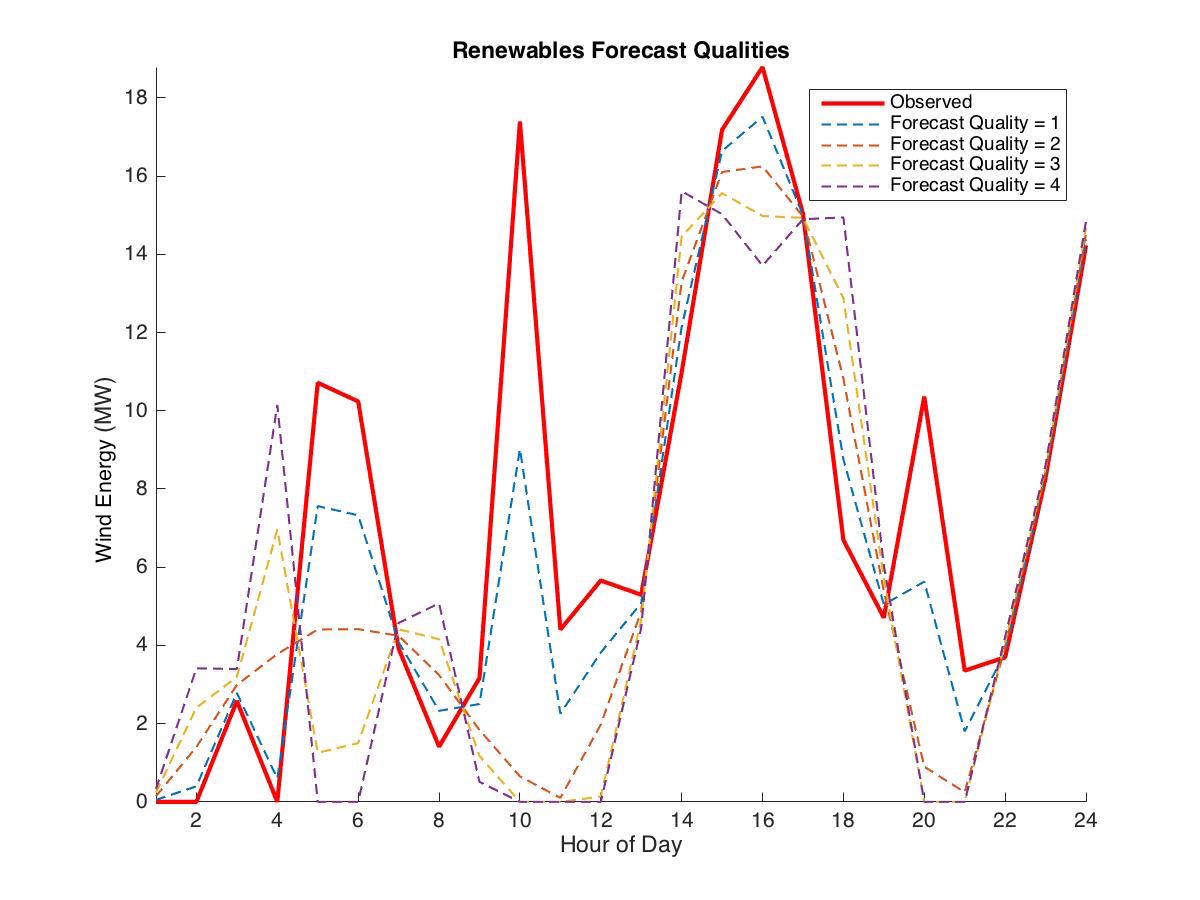}
    \caption{Forecasts of renewable energy ($E_t$)}
    \label{fig:renewables_fore}
\end{figure}

Our model below is designed in part to create complex nonstationary behaviors to test the ability of our policy to exploit forecasts while managing uncertainty. We use a hidden Markov model (\cite{durante2016Markov}) to create a very realistic model of the stochastic process describing the generation of renewable energy and make the amount of energy available from the grid a function of time. This model generates forecast errors based on an underlying crossing time distribution, the consecutive periods of time for which the observed energy produced is above or below the forecast. These errors are modeled using a two-level Markov model with two state variables that evolve on different time scales. The primary state variable, which contains all the pertinent information to approximate the current period's error distribution, evolves at every discrete point of time. The secondary state variable, also known as the crossing state of the system, contains the sign of the error and the duration of how long the sample path has been above or below the forecast. Unlike the primary state variable, this secondary state variable is only updated when forecast errors change signs. Forecast errors are then generated using a distribution selected by a second level Markov model conditioned on the crossing state of the system. A sample path of renewable energy and it's respective forecast can be viewed in figure \ref{fig:renewables}. Arrows have been added to identify crossing times. This is an example of a complex stochastic process that causes problems for stochastic lookahead models.  For example, it is very common when using the stochastic dual decomposition procedure (SDDP) to assume interstage independence, which means that $W_t$ and $W_{t+1}$ are independent, which is simply not the case in practice (\cite{shapiro2013risk} and \cite{dupavcova2002comparison}).  However, capturing this dynamic in a stochastic lookahead model is quite difficult.  Our CFA methodology, however, can easily handle these more complex stochastic models since we only need to be able to simulate the process.
We manipulate the quality of the renewable energy forecast by multiplying the forecast errors by the forecast quality, $\sigma_f$. This allows us to modify the quality of our forecast without modifying the observed stochastic process ($P_t$). Different quality forecasts for the same sample path can be seen in figure \ref{fig:renewables_fore}.
The amount of energy available from the main grid at $t$, $G_t$ is defined as:
\begin{equation}
	G_{t} = \min \bigg\{ \max \left\{ 90 - 50\sin \bigg( \frac{5\pi t}{2T} \bigg), G_{\min} \right\} G_{\max} \bigg\}  
\end{equation}
where $G_{\min}$ is the minimum energy always accessible from the grid, $G_{\max}$ is the maximum energy every accessible.

\subsection{Spot price model}
\begin{figure}[h]
    \centering
    \includegraphics[width=0.6\textwidth]{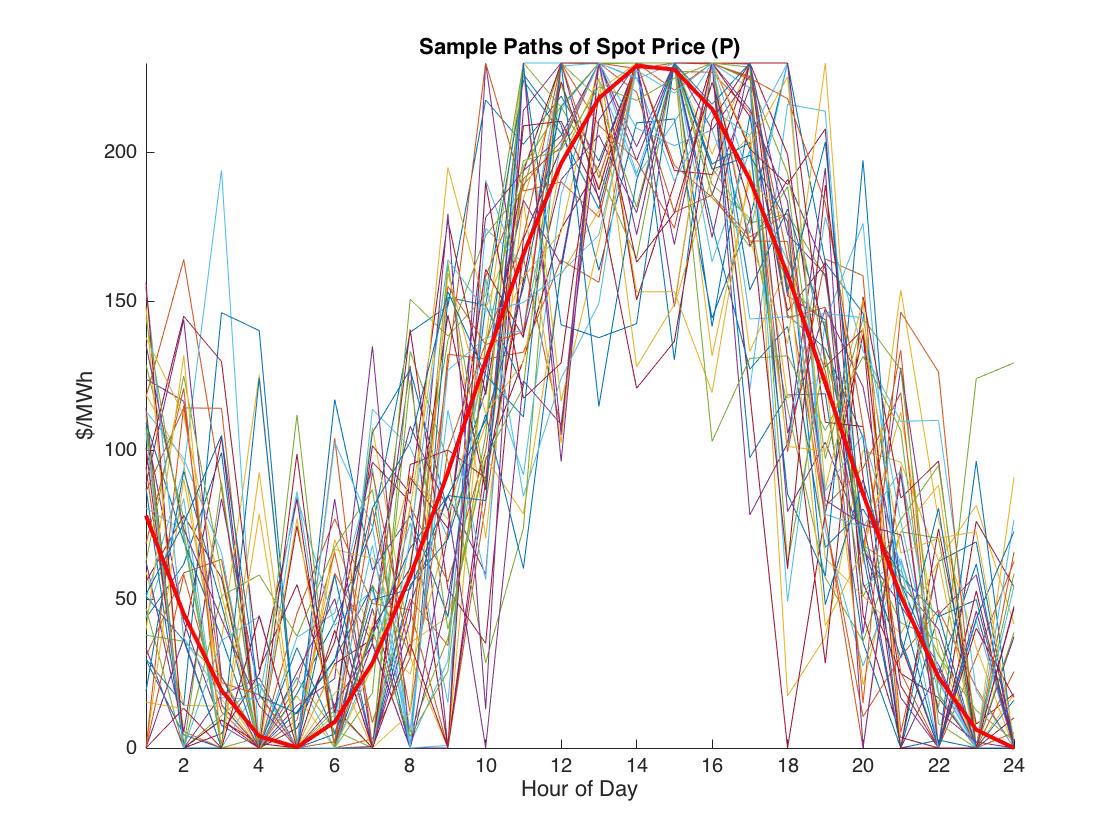}
    \caption{Sample paths of spot prices ($P_t$)}
    \label{fig:prices}
\end{figure}

The spot price ($P_t$) of electricity at time $t$ is a sinusoidal stochastic function defined as:
\begin{equation}
	P_{t} = \min \bigg\{ \max \left\{ \frac{P_{\max} + P_{\min}}{2} - (P_{\max} - P_{\min})\cdot\sin \bigg( \frac{5\pi t}{2T} \bigg) + \epsilon_t, P_{\min} \right\} P_{\max} \bigg\}  \: \:  \text{ where } \epsilon \sim \mathcal{N}(\mu_p, \sigma_p)
\end{equation}
where $P_{\min}$ is the minimum price allowed, $P_{\max}$ is the maximum price allowed, $\mu_p$ is expected value of the change in price, and $\sigma_p$ is the standard deviation of the change in price. Since spot prices occasionally go below zero $P_{\min}$ may have a negative value. This is also the price at which energy can be purchased and sold to and from the grid. Sample paths of the stochastic process $S_t$ are displayed in figure \ref{fig:prices}.

\begin{figure}[h]
    \centering
    \includegraphics[width=0.6\textwidth]{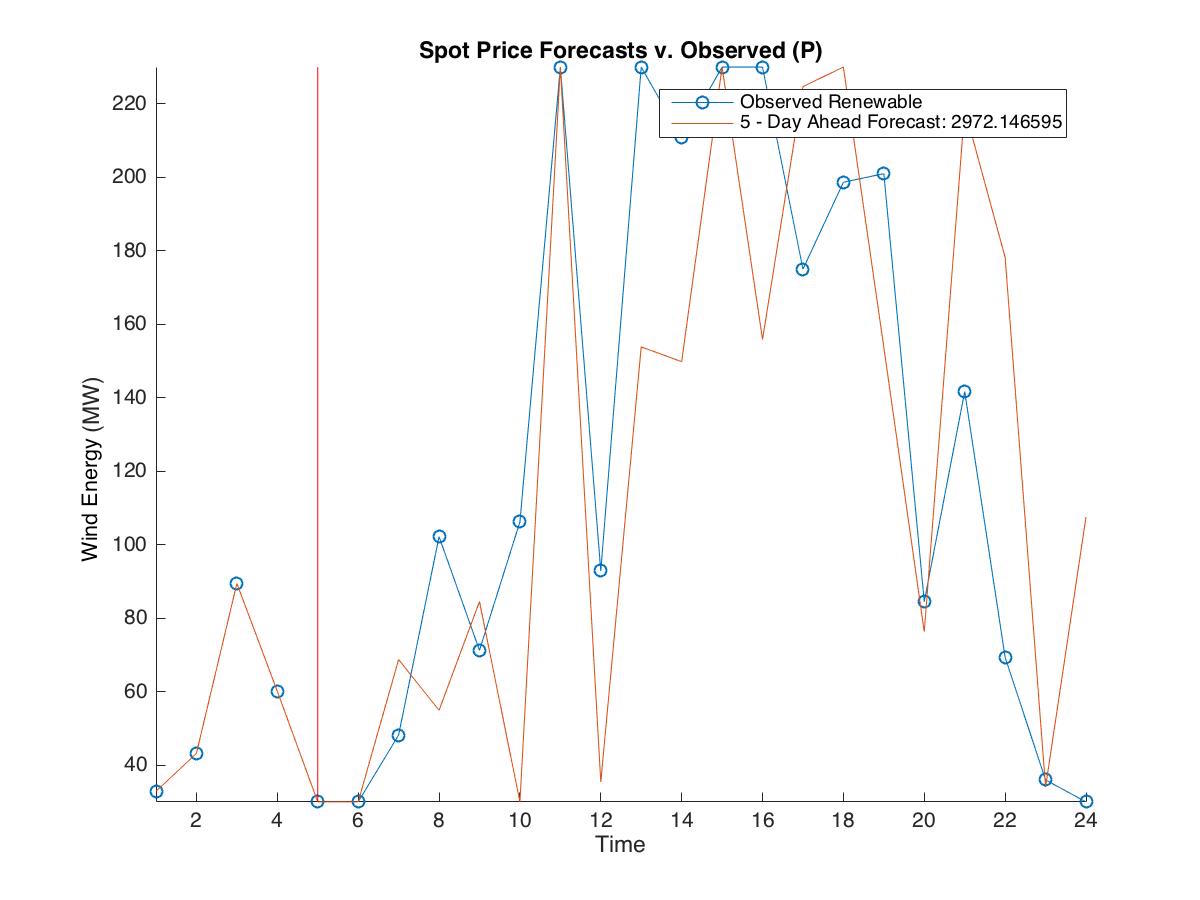}
    \caption{Spot price forecast ($P_{tt'}$) v. Observed price ($P_t$)}
    \label{fig:price_forecast}
\end{figure}

Since the price process, $P_t$, is stochastic, forecasts of $P_t$ must be generated for both the deterministic lookahead and the CFA. In our model, forecasts of spot prices are noisy observations of the observed sample path. This allows us to modify the quality of our forecast without modifying the stochastic process ($P_t$). In our simulation we generate the sample path of prices, $P_t$ where $t = 1,...,T$, first. Then we create a series of forecasts where $F^P_{tt'} = \mathbb{E}_t[P_{t'}]$ given the information available at time $t$. The process $F^P_{tt'}$ satisfies the following conditions:
\begin{enumerate}
	\item The spot price, $P_t$, is defined as
		\begin{equation}
			P_t = P_{t,t} \: \: \forall \: \: t \in [1, T]
		\end{equation}
	\item The stochastic process, $P_{tt'}$, is a Gaussian process where, 
		\begin{equation}
			P_{t-1,t'} = \min \bigg\{ \max \bigg\{\rho_t, P_{\min} \bigg\} P_{\max} \bigg\}   \: \: \: t' \geq t
		\end{equation}
\end{enumerate}
where $\rho_t \sim \mathcal{N}(P_{tt'},\sigma_f)$. We can directly control the quality of the forecast by varying $\sigma_f$, where $\sigma_f = 0$ means the forecast is perfect, while increasing $\sigma_f$ degrades the quality of the forecast. Figure \ref{fig:price_forecast} compares the forecasted price path at time $t = 5$ to the observed price path. 

\subsection{Policy Parameterizations}
	If the contribution function, transition function and constraints are linear, a deterministic lookahead policy can be constructed as a linear program if point forecasts of exogenous information are provided. For our deterministic lookahead we use the following policy  
\begin{equation} \label{energy_LA}
	X^{\text{D-LA}}_t(S_t) = \argmax_{x} \sum_{t' = t}^{t+H} c^T_{tt'} x_{tt'}
\end{equation}
{\bf subject to}
\begin{equation}\label{eqn:b_con}
		\begin{matrix}
		x^{wd}_{tt'} & + & \beta^dx^{rd}_{tt'} & + & x^{gd}_{tt'} & \leq 		& F^{D_{tt'}}, \\
		~ & ~ & x^{gd}_{tt'} & + & x^{gr}_{tt'} & \leq 	& F^{G_{tt'}}, \\
		~ & ~ & x^{rd}_{tt'} & + & x^{rg}_{tt'} & \leq 	& F^{R_{tt'}}, \\
		~ & ~ & x^{wr}_{tt'} & + & x^{gr}_{tt'} & \leq 	& R_{\max} - F^{R_{tt'}}, \\
		~ & ~ & x^{wr}_{tt'} & + & x^{wd}_{tt'} & \leq 	& F^{E_{tt'}}, \\
		~ & ~ & x^{wr}_{tt'} & + & x^{gr}_{tt'} & \leq 	& \gamma^c, \\
		~ & ~ & x^{rd}_{tt'} & + & x^{rg}_{tt'} & \leq 	& \gamma^d, \\
		\end{matrix}
\end{equation}
for $t' \in [t+1,t+H]$. We call this deterministic lookahead policy the benchmark policy, and use it to estimate the degree to which the parameterized policies are able to improve the results in the presence of uncertainty.

\begin{itemize}
	\item {\bf Capacity Constraints}: This parameterization limits the amount of energy in storage and guarantee there is capacity to purchase inexpensive energy. An upper bound constraint is easily created by multiplying the capacity of the storage device, $R_{\max}$ by the parameter $\theta_{t'-t}$. This changes the constraint 
	\begin{equation*}
		x^{wr}_{tt'} + x^{gr}_{tt'} \leq R_{\max} - F^R_{tt'}
	\end{equation*}
	to 
        \begin{equation}\label{eqn:UB}
        		x^{wr}_{tt'} + x^{gr}_{tt'} \leq R_{\max} \cdot \theta^U_{t'-t} - F^R_{tt'}
        \end{equation}
    	where $\theta_{t'} \in [0, 1]$ and $t' \in [t,t+H]$. Parameterized lower constraints are incorporated into the policy by creating the additional linear constraints
        \begin{equation}\label{eqn:LB}
        		- x_t^{rd} - x_t^{rg} + R_t \geq R_{\max} \cdot \theta^L_{t'-t}
        \end{equation}
        where $\theta^L_{t'} \in [0, 1]$ and $t' \in [t+1,t+H]$. 
	\item {\bf Lookup table forecast parameterization} - Overestimating or underestimating forecasts of renewable energy influences how aggressively a policy will store energy. We modify the forecast of renewable energy for each period of the lookahead model with a unique parameter $\theta_{\tau}$. This parameterization is a lookup table representation because there is a different $\theta$ for each lookahead period, $\tau = 0, 1, 2, ...$ The following constraints 
	\begin{equation}\label{eqn:main_con}
		x^{wr}_{tt'} + x^{wd}_{tt'} \leq F^E_{tt'}
	\end{equation}
	are changed to 
        \begin{equation}\label{eqn:UFP}
        		x^{wr}_{tt'} + x^{wd}_{tt'} \leq F^{E_{tt'}} \cdot \theta_{t' - t}.
        \end{equation}
         where $t' \in [t+1,t+H]$ and $\tau = t' -t$. If $\theta_{\tau} < 1$ the policy will be more robust and decrease the risk of running out of energy. Conversely, if $\theta_{\tau} > 1$ the policy will be more aggressive and less adamant about maintaining large energy reserves.         
	\item {\bf Constant forecast parameterization} - Instead of using a unique parameter for every period, this parameterization uses a single scalar to modify the forecast amount of renewable energy for the entire horizon. The policy constraints \eqref{eqn:main_con} are changed to 
        		\begin{equation}\label{eqn:SFP}
        			x^{wr}_{tt'} + x^{wd}_{tt'} \leq F^{E_{tt'}} \cdot \theta.
        		\end{equation}
	\item {\bf Exponential Function} - Instead of calculating a set of parameters for every period within the lookahead model we make our parameterization a function of time and two parameters. The policy constraints \eqref{eqn:main_con} are then changed to 
	\begin{equation}\label{eqn:ETD}
		x^{wr}_{tt'} + x^{gd}_{tt'} \leq F^{E_{tt'}} \cdot \theta_1 \cdot e^{\theta_2 \cdot (t'-t)}.
	\end{equation}
\end{itemize}

\section{Numerical Results}
To demonstrate the capability of the CFA and Algorithm \ref{CFA_Grad}, we test parameterizations, (\refeq{eqn:UB})-(\refeq{eqn:ETD}), of the deterministic lookahead policy defined by equation \eqref{energy_LA} on variations of the previously described energy storage problem. We provide the benchmark policy and parameterized policies the same forecasts of exogenous information. Our goal is to show parameterizing the benchmark policy and using Algorithm \ref{CFA_Grad} to determine parameter values can improve the benchmark policy's performance. We say a parameterization, $\pi(\theta)$, outperforms the nonparametric benchmark policy, if it has positive \emph{policy improvement}, $\Delta F^{\pi}(\theta)$. We define the policy improvement, $\Delta F^{\pi}(\theta)$, of parameterization $\pi(\theta)$ as
	\begin{equation}\label{improve}
		\Delta F^{\pi}(\theta) = \frac{F^{\pi(\theta)} - F^{\text{D-LA}}}{|F^{\text{D-LA}}|}
	\end{equation}
where $F^{\pi}(\theta)$ is the average profit generated by parametrization $\pi(\theta)$ and $F^{\text{D-LA}}$ is the average profit generated by the unparameterized deterministic lookahead policy described by equation \eqref{energy_LA}.        
 
        One of the most prominent advantages of the CFA is its ability to handle uncertainty without restrictions on the structure of the dynamics. By varying the forecast quality, $\sigma_f$, of the energy storage problem we demonstrate the CFA and algorithm \ref{CFA_Grad}'s ability to detect different levels of uncertainty and adapt accordingly. Table \ref{tab:SigmaF} presents the performance of each parameterization over varying forecast qualities.
        \begin{table}[h]
        \centering
        		\begin{tabular}{|l|c|c|c|c|}
\hline
&\textbf{$\sigma_f =20$}&\textbf{$\sigma_f =25$}&\textbf{$\sigma_f =30$}&\textbf{$\sigma_f =35$}\\\hline
\textbf{Constant}&13\%&13\%&16\%&17\%\\\hline
\textbf{Lookup}&20\%&22\%&26\%&25\%\\\hline
\textbf{Expo}&14\%&22\%&26\%&26\%\\\hline
\textbf{Capacity Con.}&0.00\%&0.00\%&0.00\%&0.00\%\\\hline
\end{tabular}

        	\caption{This table displays the percentage improvement obtained by parameterized policies relative to the deterministic benchmark for varying forecast qualities, $\sigma_f$. These values are calculated using $500$ simulations.}
         \label{tab:SigmaF}
         \end{table}

\begin{figure}[h]
                  \centering
                  \includegraphics[width=.75\linewidth]{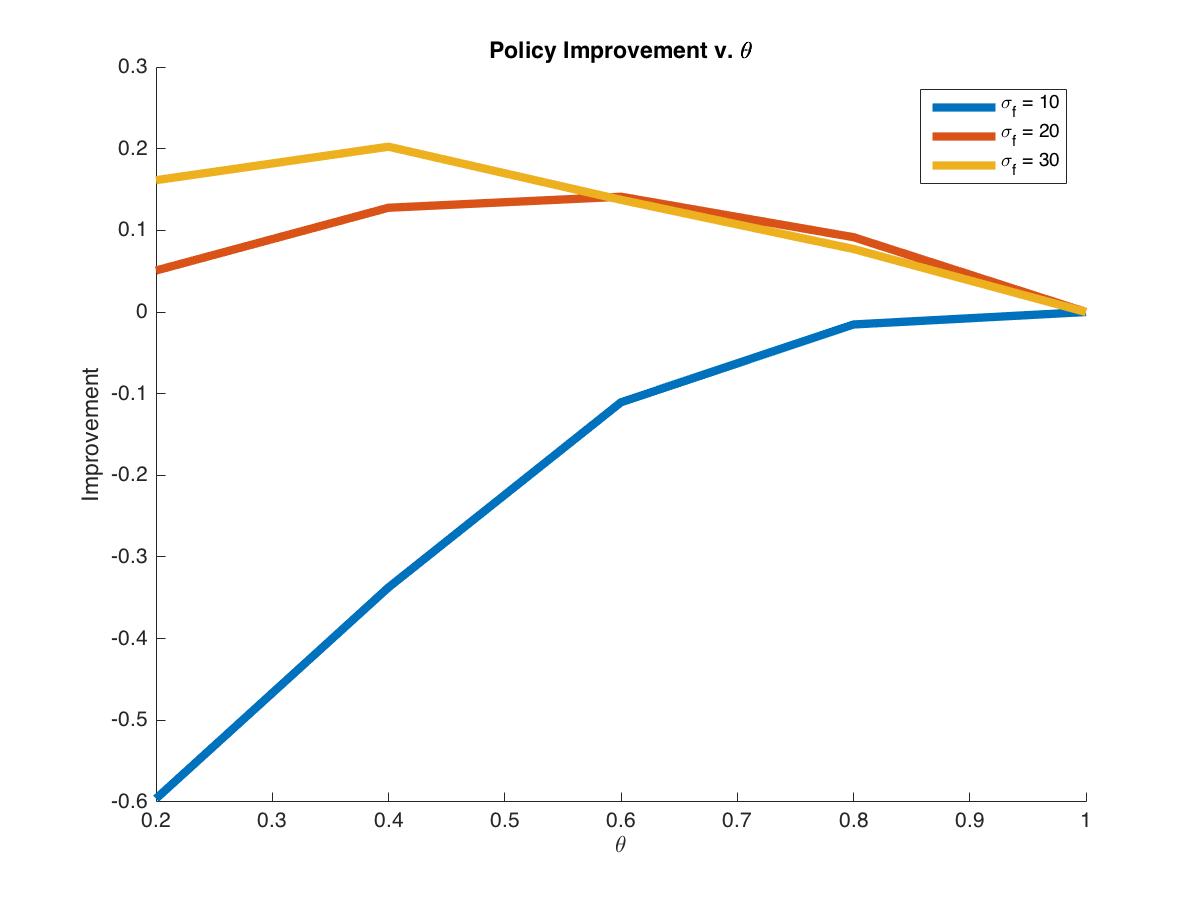}
                  \caption{Policy Improvement over deterministic benchmark v. $\theta$ for constant parameterizations}
                  \label{fig:Sigma_single}	
         \end{figure}

Given a perfect forecast, $\sigma_f = 0$, the benchmark policy, a deterministic lookahead, is the optimal policy. As the uncertainty and forecast error increases the benchmark policy's performance deteriorates and the average profit generated decreases y since it is unable to deal with uncertainty.  The average profit of the parameterizations also deteriorate as forecast error increase, but does so at a slower rate than the benchmark policy. Although the added noise to the forecast makes the problem more difficult, the parameterized policy is able to adapt and perform better than the standard deterministic lookahead policy. This explains the positive relationship between the the Constant, Lookup table, and Exponential parameterizations improvements and forecast quality. As the forecast error, $\sigma_f$, increases the less trustworthy the forecast becomes. To account for this increased uncertainty the policies further underestimates the forecast to limit the risk of paying penalties for not satisfying demand. This phenomena can be seen in figure \ref{fig:Sigma_single}, a plot of the relationship of $\theta$ for the constant parameterization and forecast qualities.                    
               
%
%

\subsection{Factoring the Forecast}
\begin{figure}
        \centering
        \begin{subfigure}{.5\textwidth}
          \centering
         \includegraphics[width=1\linewidth]{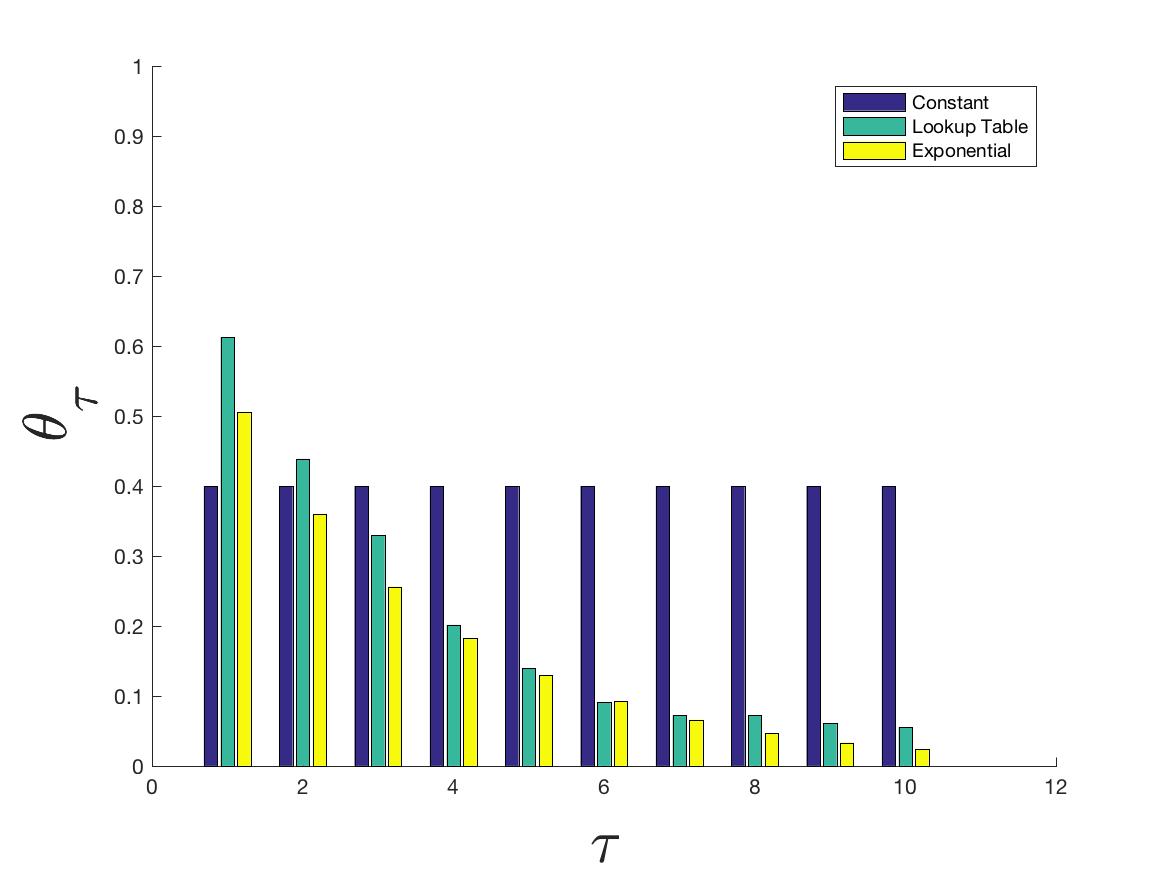}
         \caption{$\theta_\tau$ for $\sigma_f = 35$ and $\gamma = R_{\max}$}
          \label{fig:ThetaBar_gamma}
        \end{subfigure}%
        \begin{subfigure}{.5\textwidth}
          \centering
          \includegraphics[width=1\linewidth]{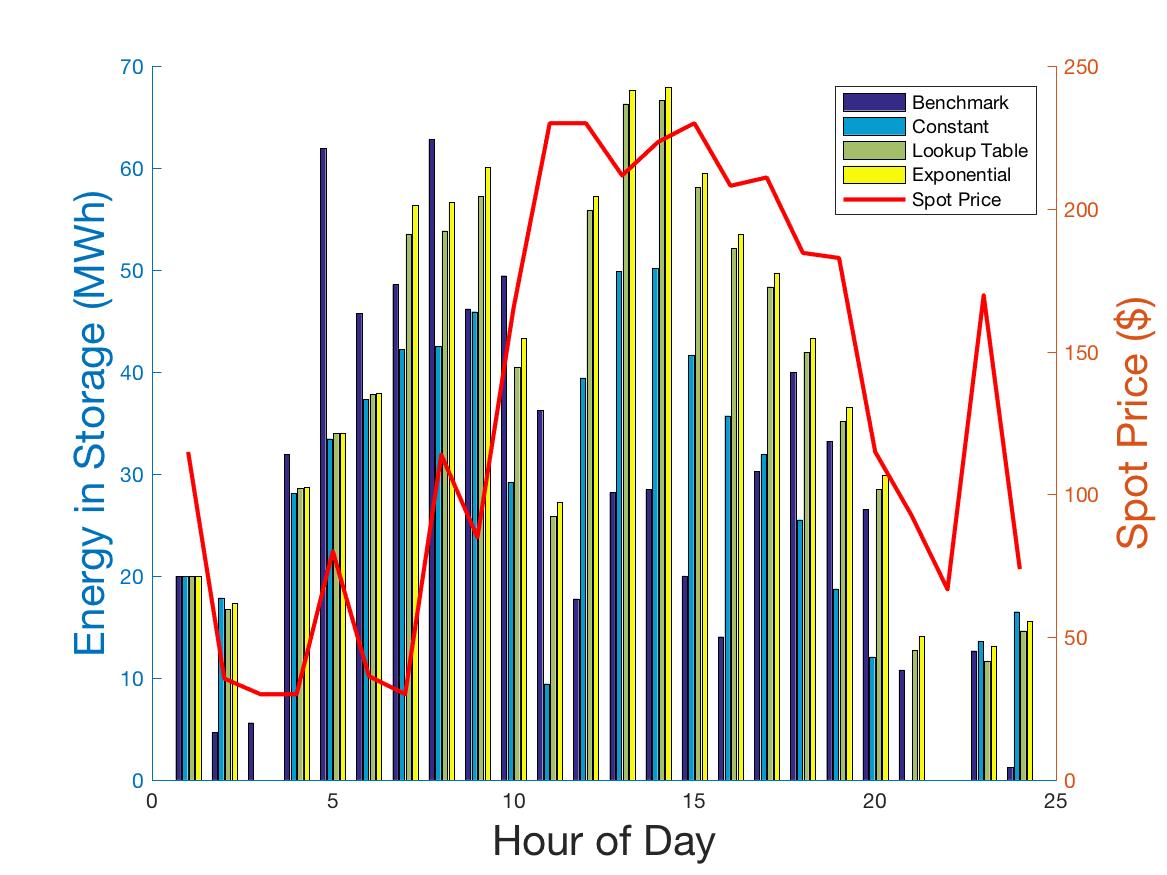}
          \caption{$R_t$ for $\sigma_f = 35$ and $\gamma = .R_{\max}$}
          \label{fig:Rbar1}
        \end{subfigure}
        \caption{Figure \ref{fig:ThetaBar_gamma} compares the $\theta_\tau$ values for the \emph{Constant, Lookup Table}, and \emph{Exponential} parameterizations when $\sigma_f = 35$ and $\gamma = \cdot R_{\max}$. Figure \ref{fig:Rbar1} compares the storage levels, $R_t$ of the different parameterizations over $t \in [1,24]$ for the same conditions.}
        \label{fig:ThetaBar_sigmaf}
        \end{figure}

All policies prefer to underestimate estimate the future renewable levels by setting $\theta_{\tau} < 1$ for all $\tau \in [0,H]$. Figures \ref{fig:ThetaBar_sigmaf} and \ref{fig:ThetaBar_gamma} shows how $\theta_\tau$ for the parameterizations described by \eqref{eqn:UFP} - \eqref{eqn:ETD} behave as functions of $\tau$.
        
	
        	\begin{figure}
        \centering
        \begin{subfigure}{.5\textwidth}
          \centering
         \includegraphics[width=1\linewidth]{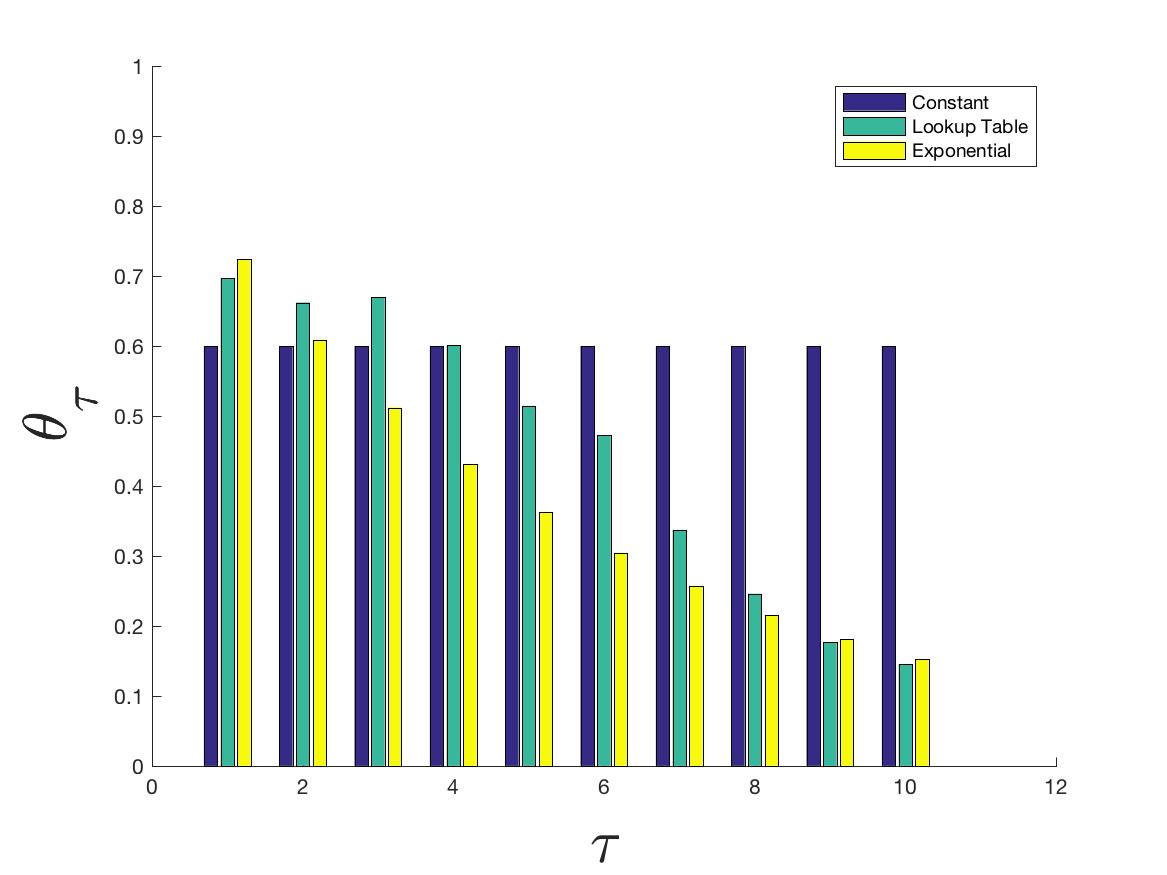}
         \caption{$\theta_\tau$ for $\sigma_f = 20$ and $\gamma = R_{\max}$}
          \label{fig:ThetaBar_gamma2}
        \end{subfigure}%
        \begin{subfigure}{.5\textwidth}
          \centering
          \includegraphics[width=1\linewidth]{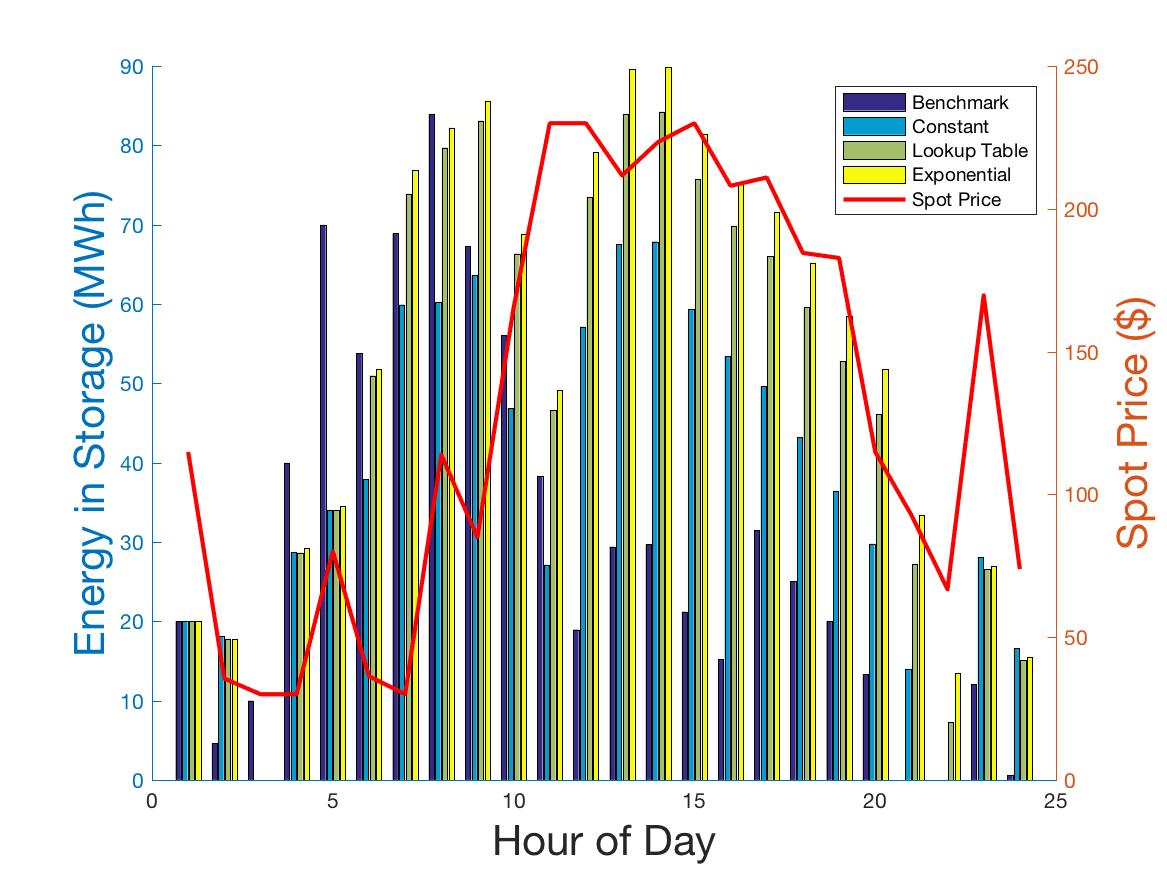}
          \caption{$R_t$ for $\sigma_f = 20$ and $\gamma = R_{\max}$}
          \label{fig:Rbar2}
        \end{subfigure}
        \caption{Figure \ref{fig:ThetaBar_gamma2} compares the $\theta_\tau$ values for the \emph{Constant, Lookup Table}, and \emph{Exponential} parameterizations when $\sigma_f = 30$ and $\gamma = R_{\max}$. Figure \ref{fig:Rbar2} compares the storage levels, $R_t$ of the different parameterizations over $t \in [1,24]$ for the same conditions.}
        \label{fig:ThetaBar_gamma}
        \end{figure}

Notice how $\theta_\tau$ decreases for each subsequent lookahead period for the lookup table and exponential adjustment functions. This a consequence of the diminishing marginal improvement for each additional period in the lookahead model. As seen in figure \ref{fig:Sigma_single}, as the forecast error, $\sigma_f$, increases, $\theta_\tau$ deteriorates. This implies the algorithm recognizes that as the forecast error increases the forecast is less reliable. The policy determines that it is better to just expect no renewable energy than to depend the forecast.
        
\subsection{Capacity Constraints}
\begin{figure}[h]
        \centering
        		\includegraphics[width=1\linewidth]{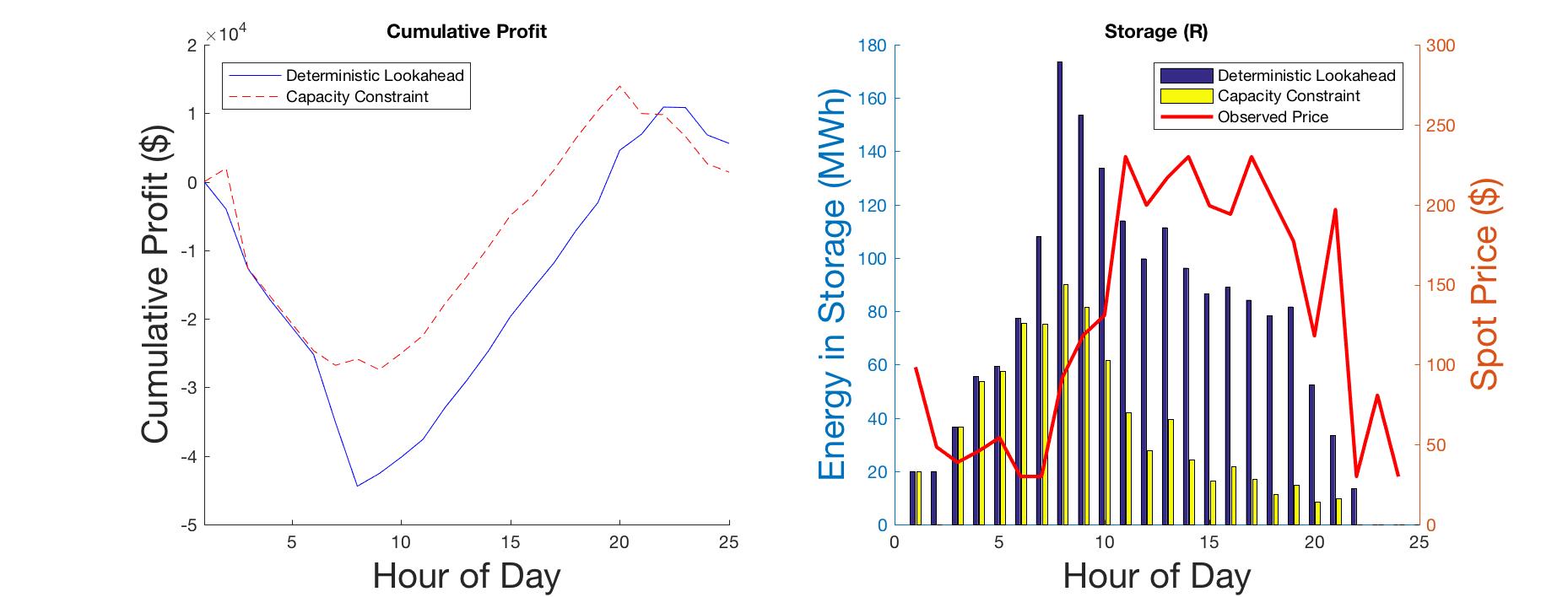}
        	\caption{Capacity constraint parameterization sample path where $\theta = [.3,0]$}
         \label{fig:CC_samp1}
         \end{figure}

 The capacity constraint parameterization, described by equations \eqref{eqn:UB} and \eqref{eqn:LB}, was the only parameterization that did not generate positive improvement in the provided problem settings. Setting an upper limit on the lookahead model's storage decreases the amount of energy placed into storage during the current state. This maintains lower storage levels than the benchmark policy and limits the purchased energy from the grid for storage. However, this also limits the parameterized policy's opportunity to sell excess energy to the grid for profit. This can be seen in figure \ref{fig:CC_samp1}.

\begin{figure}[h]
        \centering
        		\includegraphics[width=1\linewidth]{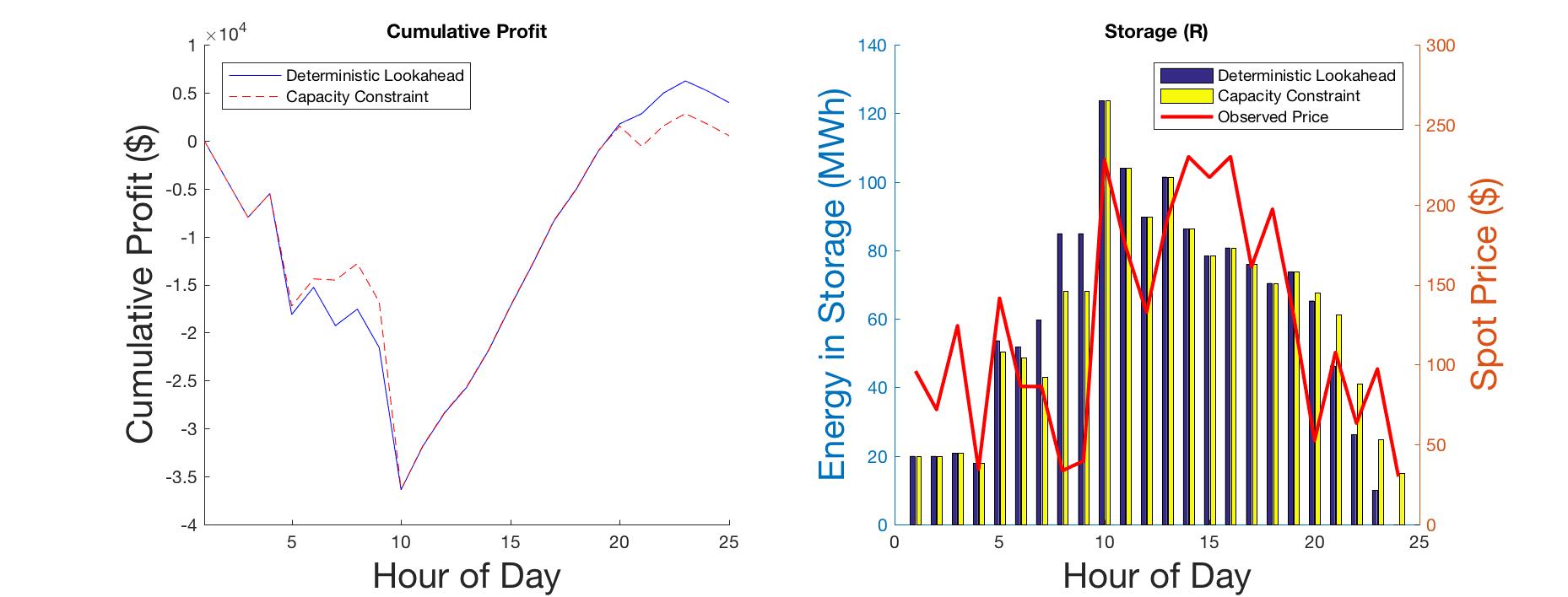}
        	\caption{Capacity constraint parameterization sample path where $\theta = [1,.05]$}
         \label{fig:CC_samp2}
         \end{figure}

Notice how the parameterized policy's cumulative profit is greater than the benchmark's until $t = 20$ in figure \ref{fig:CC_samp1}. The parameterized policy achieves this by maintaining lower storage associated costs. However, as the simulation approaches $t = T$ the benchmark policy begins to sell off excess storage. Since the parameterized policy's storage is constantly lower than the benchmark it misses the additional returns. Setting a lower limit has the reverse effect on storage. This can be seen in figure \ref{fig:CC_samp2}.

By requiring a certain amount of energy in the lookahead model's storage device the policy is unable to sell as much excess energy to the grid as the benchmark policy. This limits the policy's ability to generate revenue. Algorithm \ref{CFA_Grad} seemed to recognize these problems and did not limit the lookahead model's capacity constraints much. Though it could not improve the lookahead policy by modifying the capacity constraints, it still identified the optimal $\theta^* = [1,0]$ for the parameterization form. 

\section{Conclusion}
This work builds upon a long history of using deterministic optimization models to solve sequential stochastic problems. Unlike other deterministic methods, our class of methods, {\it parametric cost function approximations}, parametrically modify deterministic approximations to account for problem uncertainty. Our particular use of modified linear programs and the  Gradient Algorithm in Algorithm \ref{CFA_Grad} are fundamentally new. Our method allows us to exploit structural properties of the problem and compensate for uncertainty simultaneously. We have demonstrated this class of policies in the context of a very rich class of energy storage problems. For our numerical work we selected an energy storage problem that is relatively small to simplify the extensive computational work. However, our methodology is scalable to any problem setting which is currently being solved using a deterministic model.

This new class of policies offers a new breadth of research possibilities such as identifying other appropriate problem classes and policy structures.  We also recognize that gradient-based search mechanism are not always possible. Therefore developing derivative-free stochastic search methods for tuning CFAs is another potential area of future work, as well as designing methods to do adaptive search in an online setting.

\newpage 
\bibliographystyle{agsm}
\bibliography{library}
\end{document}